\documentclass[11pt,a4paper]{article}

\usepackage{geometry}
\usepackage[latin1]{inputenc}
\usepackage{amsmath,amsfonts,amssymb,amsthm}
\usepackage{mathrsfs}
\usepackage{bm,bbm}
\usepackage{booktabs}
\usepackage{enumitem}

\usepackage{lscape}

\usepackage[colorlinks=true,linkcolor=red,citecolor=blue]{hyperref}

\usepackage{natbib}
\usepackage{graphicx,epsfig}
\graphicspath{ {./Images/} }
\usepackage{caption}
\usepackage{subcaption}

\usepackage{xcolor}
\definecolor{electricviolet}{rgb}{0.56, 0.0, 1.0}

\numberwithin{equation}{section}

\newtheorem{theorem}{Theorem}
\newtheorem{proposition}{Proposition}
\newtheorem{corollary}{Corollary}

\theoremstyle{definition}
\newtheorem{remark}{Remark}

\DeclareMathOperator*{\argmin}{arg\,min}
\newcommand{\ind}{\mathbbm{1}}

\begin{document}

\title{Extreme expectile estimation for short-tailed data, with an application to market risk assessment}

\author{Abdelaati Daouia$^{a}$, Simone A. Padoan$^{b}$ \& Gilles Stupfler$^{c}$}
\date{$^{a}$ {\small Toulouse School of Economics, University of Toulouse Capitole, France} \\ $^{b}$ {\small Department of Decision Sciences, Bocconi University, via Roentgen 1, 20136 Milano, Italy} \\ $^{c}$ {\small Univ Angers, CNRS, LAREMA, SFR MATHSTIC, F-49000 Angers, France}}

\maketitle 

%
\begin{abstract} 
The use of expectiles in risk management has recently gathered remarkable momentum due to their excellent axiomatic and probabilistic properties. In particular, the class of elicitable law-invariant coherent risk measures only consists of expectiles. While the theory of expectile estimation at central levels is substantial, tail estimation at extreme levels has so far only been considered when the tail of the underlying distribution is heavy. This article is the first work to handle the short-tailed setting where the loss ({\it e.g.} negative log-returns) distribution of interest is bounded to the right and the corresponding extreme value index is negative. We derive an asymptotic expansion of tail expectiles in this challenging context under a general second-order extreme value condition, which allows to come up with two semiparametric estimators of extreme expectiles, and with their asymptotic properties in a general model of strictly stationary but weakly dependent observations. A simulation study and a real data analysis from a forecasting perspective are performed to verify and compare the proposed competing estimation procedures.
\end{abstract}



\section{Introduction}
\label{sec:intro}

The class of expectiles, introduced by~\cite{newpow1987}, defines 
useful descriptors $\xi_{\tau}$ of the higher ($\tau\ge \frac{1}{2}$) and lower ($\tau\le \frac{1}{2}$) regions of the distribution of a random variable~$X$ through the
asymmetric least squares (ALS) minimization problem 
\begin{equation}
\label{eqn:expectile}
\xi_{\tau}= \argmin_{\theta\in\mathbb{R}} \mathbb{E} \big[ \eta_{\tau}(X-\theta) - \eta_{\tau}(X) \big] ,
\end{equation}
where $\eta_{\tau}(x) = | \tau-\ind\{ x\leq 0 \}|\,x^2$, with $\tau\in (0,1)$ and $\ind\{ \cdot \}$ being the indicator function.
Expectiles are well-defined, finite and uniquely determined as soon the first moment of $X$ is finite. They 
generalize the mean 
$\xi_{1/2}=\mathbb{E}(X)$
in the same way quantiles generalize the median, thus defining an ALS analog to quantiles. Indeed,
 \cite{koebas1978} showed that the $\tau$th quantile $q_{\tau}$ of $X$ solves the asymmetric $L^1$ minimization problem 
\[
q_{\tau}\in \argmin_{\theta\in\mathbb{R}} \mathbb{E} ( \varrho_{\tau}(X-\theta) - \varrho_{\tau}(X) ) ,
\]
where $\varrho_{\tau}(x) = | \tau-\ind\{ x\leq 0 \}|\,|x|$.  
Expectiles have received renewed attention for their ability to quantify tail risk 
at least since the contribution of~\cite{tay2008}. They depend on the tail realizations of $X$ and their probability, while quantiles only depend on the frequency of tail realizations, see~\cite{kuayehhsu2009}. Most importantly, \cite{zie2016} showed that expectiles are the sole coherent law-invariant measure of risk which is also elicitable in the sense of~\cite{gne2011}, meaning that they abide by the intuitive diversification principle~\citep{belklamulgia2014} and that their prediction can be performed through a straightforward principled backtesting methodology.
These merits have motivated the development of procedures for expectile estimation and inference over the last decade. A key, but difficult, question in any risk management setup is the estimation of the expectile $\xi_{\tau}$ at extreme levels, which grow to 1 as the sample size increases. This question was first tackled in~\cite{daogirstu2018,daogirstu2020} under the assumption that the underlying distribution is heavy-tailed, that is, its distribution function tends to 1 algebraically fast. The latest developments under this assumption have focused on, among others, bias reduction~\citep{girstuauc2022}, accurate inference~\citep{padstu2022}, and handling more complex data in regression~\citep{girstuauc2021,girstuauc2022} or time series~\citep{davpadstu2022} setups. 

The problem of estimating extreme expectiles outside of the set of heavy-tailed models is substantially more complicated from a statistical standpoint. The contribution of the present paper is precisely to build and analyze semiparametric extreme expectile estimators in the challenging short-tailed model, in which the extreme value index (EVI) of the underlying distributions is known to be negative. This requires employing a dedicated extrapolation relationship for population extreme expectiles. Only~\cite{maonghu2015} have initiated such a study at the population level when~$X$ belongs to the domain of attraction of a Generalized Extreme Value distribution (GEV). Differently to~\cite{maonghu2015}, we work in the general semiparametric Generalized Pareto (GP) setting through a standard second-order condition, which makes it possible to 
 derive an asymptotic expansion of extreme expectiles without resorting to an unnecessary restriction about the link between the 
 EVI and second-order parameter that featured in~\cite{maonghu2015}. Based on this asymptotic expansion, we present and study two different extreme value estimators of tail expectiles. The first one builds upon the Least Asymmetrically Weighted Squares (LAWS) estimator of expectiles, namely the empirical counterpart of $\xi_{\tau}$ in~\eqref{eqn:expectile}, obtained at intermediate levels $\tau=\tau_n\to 1$ with $n(1-\tau_n)\to \infty$ as the sample size $n\to\infty$. 
The short-tail model assumption allows then to come up with an expectile estimator extrapolated to the far tail at arbitrarily extreme levels $\tau=1-p_n$ such that $(1-\tau_n)/p_n\to \infty$ as $n\to\infty$, in a semiparametric way reminiscent of how extreme quantiles are fitted in Section~4.3 of~\cite{haafer2006}.
The second extrapolating estimator directly relies on the asymptotic expansion of $\xi_{\tau}$ that involves 
its quantile analog $q_{\tau}$, the endpoint $q_1\equiv\xi_1$ and the EVI, by plugging in the GP quantile-based estimators of these tail quantities.
Our estimation theory is valid in a general setting of strictly stationary and weakly dependent 
data satisfying reasonable mixing and tail dependence conditions.  
We explore various theoretical and practical features of extreme expectile estimation 
in this setting,
and explain why this problem is statistically more difficult than extreme quantile estimation. In particular, 
an extreme expectile $\xi_{\tau}$ is intrinsically less spread than its quantile analog $q_{\tau}$, even at asymmetry levels $\tau\approx 1$ where it remains much closer to the center of the distribution than~$q_{\tau}$. 
Consequently, any semiparametric 
procedure for extreme expectile estimation should be expected to suffer at least from a worse bias than for extreme quantile estimation.

Our focus on the problem of estimating extreme expectiles for bounded distributions
is motivated by the perhaps somewhat surprising finding that weekly returns of equities, used in applications to circumvent the non-synchronicity of daily data, may have short-tailed distributions.
%
\begin{figure}[t!] 
\centering
    \includegraphics[width=14cm, height=10cm]{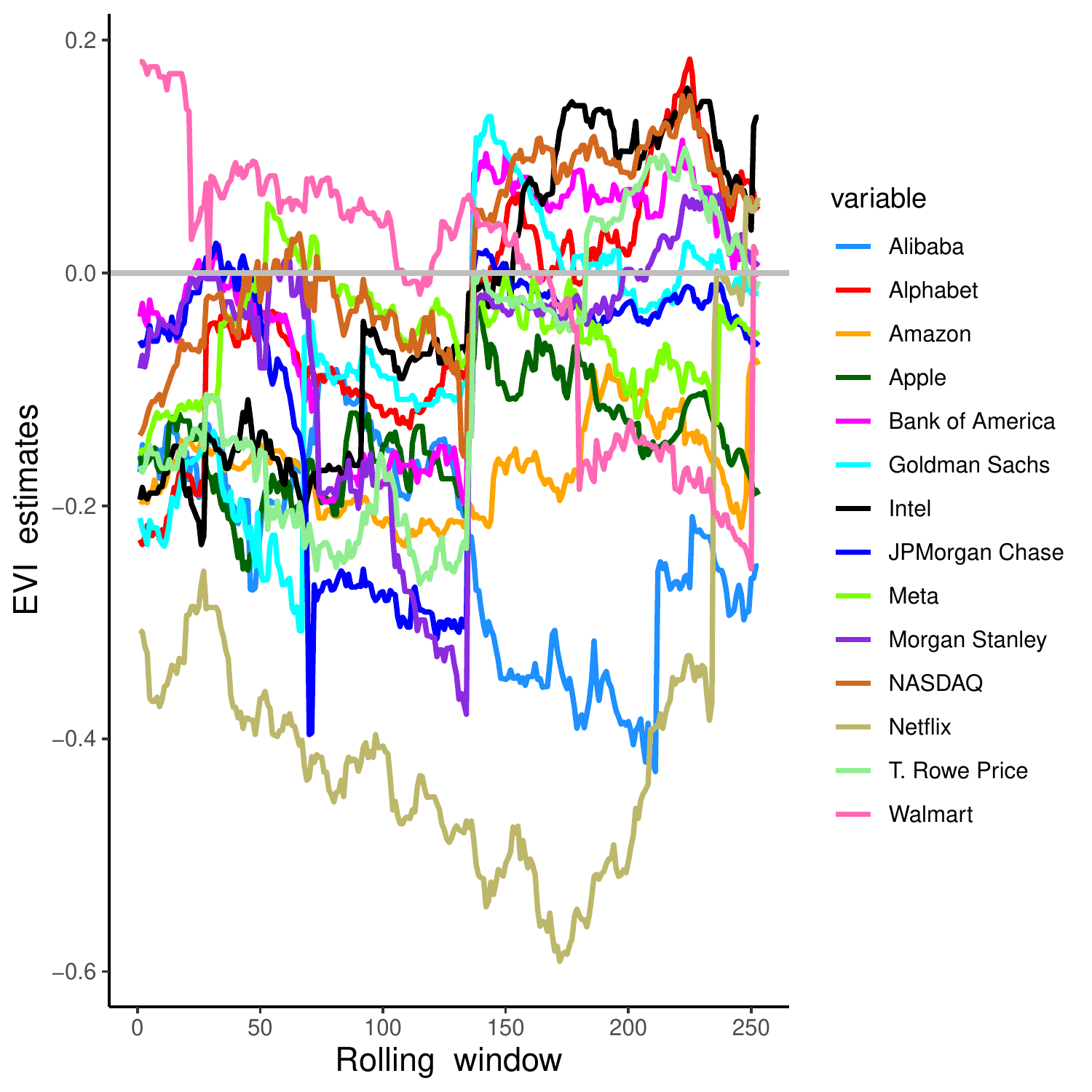} 
\caption{Maximum Likelihood estimates of the 
extreme value index over the resulting 253 successive rolling windows of~$150$ stationary data, obtained from 14 time series of weekly logarithmic (loss) returns between 21st September 2014, and 12th June 2022.}
\label{fig:data:EVI}
\end{figure}
%
This is illustrated in Figure~\ref{fig:data:EVI}  
for~14 major companies and financial institutions, where the data consists of the loss returns ({\it i.e.}~negative log-returns) on their weekly equity price from 21st September 2014 to 12th June 2022, corresponding to 403 trading weeks. The representative price is constructed by averaging daily closing prices within the corresponding week. 
The nature of the upper tail of these loss returns is reflected by the 
EVI of their distribution
whose negative, zero or positive values indicate respectively a distribution with short, light or Pareto-type tail. None of these three scenarios can be excluded in practice for these 14 data examples, where the EVI
is estimated on successive rolling windows of length $n=150$ using the GP 
distribution
fitted to exceedances over a high threshold by means of the 
Maximum Likelihood (ML) method, with the optimal threshold being chosen by
the path stability procedure as described below in Section~\ref{sec:data}. 
It is therefore important 
to construct an appropriate and fully data-driven estimation procedure for the challenging scenario of short-tailed data. 
This problem also appears naturally in production econometrics when analyzing the productivity of firms~\citep{Kokic1997}.
All our methods and data have been incorporated into the {\tt R} package {\tt ExtremeRisks}.  

In Section~\ref{sec:main}, we explain in detail the short tail distributional assumption on $X$, state our asymptotic expansion linking extreme expectiles and quantiles, construct our two classes of extreme expectile estimators and study their asymptotic properties. A simulation study examines their finite-sample performance in Section~\ref{sec:fin}, and a 
time series of 
Bitcoin data is analyzed in Section~\ref{sec:data}. 
The online Supplementary Material contains all the proofs in Section A and further simulation results in Section B.

\section{Main results}
\label{sec:main}

\subsection{Connection between extreme expectiles and quantiles}
\label{sec:main:connect}

Let $F:x\mapsto \mathbb{P}(X\leq x)$ be the distribution function of the random variable 
of interest $X$
and $\overline{F}=1-F$ be 
its
survival function. Define the associated quantile function by $q_{\tau}=\inf\{ x \in \mathbb{R} \, | \, F(x)\geq \tau\}$ and the tail quantile function $U$ by $U(s)=q_{1-s^{-1}}$, $s>1$. Differently from existing literature on extreme expectile estimation, we focus on the case when the distribution of $X$ is short-tailed,
or equivalently, when 
its EVI $\gamma$
is negative. According to Theorem~1.1.6 on p.10 of~\cite{haafer2006}, this corresponds to assuming that there is a positive function $a$ such that
\[
\forall z>0, \ \lim_{s \to \infty} \frac{U(sz)-U(s)}{a(s)} = \frac{z^{\gamma}-1}{\gamma}, \mbox{ with } \gamma<0. 
\]
This assumption can be informally rewritten as
\begin{equation}
\label{eqn:approx}
\forall z>0, \ U(sz) \approx U(s) + a(s) \frac{z^{\gamma}-1}{\gamma} \mbox{ when } s \mbox{ is large.}
\end{equation}
This means that 
extreme values of $X$ at the far tail (represented by $U(sz)$) can be achieved by extrapolating
in-sample large values (represented by $U(s)$) if the scale function $a(s)$ and the shape parameter $\gamma$ can be consistently estimated. 
The theory of the resulting extreme value estimators is usually developed under the following second-order refinement of the short-tailed model assumption above, which will be our main condition throughout~\citep[see][Equation~(2.3.13) on~p.45]{haafer2006}:
\vskip1ex
\noindent
Condition $\mathcal{C}_2( \gamma, a, \rho, A) \ $ There exist $\gamma<0$, $\rho \leq 0$, a positive function $a(\cdot)$ and a measurable function $A(\cdot)$ having constant sign and converging to $0$ at infinity such that, for all $z>0$, 
\[
\lim_{s \to \infty} \frac{1}{A(s)} \left( \frac{U(sz)-U(s)}{a(s)}-\frac{z^{\gamma}-1}{\gamma} \right) = \int_1^z v^{\gamma-1} \left( \int_1^v u^{\rho-1} \mathrm{d}u \right) \mathrm{d}v.
\]

\vskip1ex
\noindent
This condition enables one to control the bias incurred by using 
the approximation~\eqref{eqn:approx} and 
represented by the function $A$. Under this condition, the right endpoint $x^{\star}=\sup\{ x\in \mathbb{R} \, | \, F(x)<1 \}$ of $X$ is necessarily finite~\citep[see][Theorem~1.2.1 on~p.19]{haafer2006}. 
This justifies calling this model a short-tailed (or bounded) model.

Suppose now that $\mathbb{E}|\min(X,0)|<\infty$ and that condition $\mathcal{C}_2( \gamma, a, \rho, A)$ is satisfied, so that $\mathbb{E}|X|<\infty$ and expectiles of $X$ are well-defined and finite. 
First, we
motivate an asymptotic expansion of extreme expectiles that will be instrumental in our subsequent theory of extreme expectile estimation. Recall that the $\tau$th expectile $\xi_{\tau}$ satisfies 
\begin{equation}
\label{eqn:expectiledefinition}
\xi_{\tau}-\mathbb{E}(X)=\frac{2\tau-1}{1-\tau}\mathbb{E}((X-\xi_{\tau})\ind\{X>\xi_{\tau}\}),
\end{equation}
see~Equation~(12) in~\cite{belklamulgia2014}. Writing $\mathbb{E}((X-x)\ind\{X>x\})$ as an integral of the quantiles of $X$ above $x$ and using condition $\mathcal{C}_2( \gamma, a, \rho, A)$ justifies the approximation
\[
\mathbb{E}((X-\xi_{\tau})\ind\{X>\xi_{\tau}\}) \approx \frac{\overline{F}(\xi_{\tau}) a(1/\overline{F}(\xi_{\tau}))}{1-\gamma} \ \mbox{ as } \tau\uparrow 1,
\]
and therefore 
\begin{equation}
\label{eqn:asyeqF}
\lim_{\tau\uparrow 1} \frac{a(1/\overline{F}(\xi_{\tau})) \overline{F}(\xi_{\tau})}{1-\tau} = (x^{\star}-\mathbb{E}(X))(1-\gamma).
\end{equation}
%
The convergence $a(s)/(x^{\star} - U(s))\to -\gamma$ as $s\to \infty$~\citep[see][Lemma~1.2.9 on~p.22]{haafer2006} then suggests
\begin{equation}
\label{eqn:asyeqF2}
\lim_{\tau\uparrow 1} \frac{(x^{\star} - \xi_{\tau}) \overline{F}(\xi_{\tau})}{1-\tau} = (x^{\star}-\mathbb{E}(X))(1-\gamma^{-1}).
\end{equation}
The approximations $\overline{F}(\xi_{\tau})/(1-\tau)\approx \overline{F}(\xi_{\tau})/\overline{F}(q_{\tau})\approx (x^{\star} - \xi_{\tau})^{-1/\gamma}/(x^{\star} - q_{\tau})^{-1/\gamma}$ motivated by the regular variation property of $x\mapsto \overline{F}(x^{\star}-1/x)$~\cite[see][Theorem~1.2.1.2 on~p.19]{haafer2006} finally entail
%
\begin{equation}
\label{eqn:asyeq}
\lim_{\tau\uparrow 1} \frac{x^{\star}-\xi_{\tau}}{(x^{\star}-q_{\tau})^{1/(1-\gamma)}} = [(x^{\star}-\mathbb{E}(X))(1-\gamma^{-1})]^{-\gamma/(1-\gamma)}.
\end{equation}
Consequently, extreme expectiles can be extrapolated from their quantile analogs in conjunction with endpoint and 
EVI estimation.
Analyzing the asymptotic properties of 
the estimators built in this way will require quantifying the difference between the ratio $(x^{\star}-\xi_{\tau})/(x^{\star}-q_{\tau})^{1/(1-\gamma)}$ and its limit in~\eqref{eqn:asyeq}. This is the focus of our first main result below.
\begin{proposition}
\label{prop:asyexpansion}
Suppose that $\mathbb{E}|\min(X,0)|<\infty$ and condition $\mathcal{C}_2( \gamma, a,\rho, A)$ holds with $\rho<0$, and let $x^{\star}$ be the finite right endpoint of $F$. Then 
\begin{align*}
x^{\star} - \xi_{\tau} &= 
 \frac{(x^{\star} - q_{\tau})^{\frac{1}{1-\gamma}}}{[(x^{\star}-\mathbb{E}(X))(1-\gamma^{-1})]^{\frac{\gamma}{1-\gamma}}} 
    \Bigg( 1 - \left(\frac{x^{\star}-q_{\tau}}{(x^{\star}-\mathbb{E}(X))(1-\gamma^{-1})}\right)^{\frac{1}{1-\gamma}}(1+\operatorname{o}(1)) \\
    &\qquad +\frac{\gamma[(x^{\star}-\mathbb{E}(X))(1-\gamma^{-1})]^{-\frac{\rho}{1-\gamma}}}{\rho(\gamma+\rho)(1-\gamma-\rho)} A\left(\frac{(x^{\star} - q_{\tau})^{\frac{1}{1-\gamma}}}{1-\tau}\right)(1+\operatorname{o}(1)) \Bigg)
\end{align*}
as $\tau\uparrow 1$. In particular 
\begin{multline*}
    x^{\star} - \xi_{\tau} = \frac{(x^{\star} - q_{\tau})^{\frac{1}{1-\gamma}}}{ [(x^{\star}-\mathbb{E}(X))(1-\gamma^{-1})]^{\frac{\gamma}{1-\gamma}}} 
    \left( 1 + \operatorname{O}((1-\tau)^{-\frac{\gamma}{1-\gamma}}) + \operatorname{O}(|A((1-\tau)^{-\frac{1}{1-\gamma}})|) \right).
\end{multline*}
\end{proposition}
The additional condition $\rho<0$ in Proposition~\ref{prop:asyexpansion} is very mild and satisfied in all standard short-tailed models, see~\citet[][Table~2.2 on~p.68]{beigoesegteu2004}.
\begin{remark}[An equivalent asymptotic expansion]
\label{rmk:asyexpansion}
An equivalent version of the asymptotic expansion in Proposition~\ref{prop:asyexpansion}, where the error terms rely solely on the expectile $\xi_{\tau}$, is 
\begin{multline*}
    x^{\star} - \xi_{\tau} = [(x^{\star}-\mathbb{E}(X))(1-\gamma^{-1})]^{-\gamma/(1-\gamma)} (x^{\star} - q_{\tau})^{1/(1-\gamma)} \\
    \times \left( 1 - \frac{x^{\star}-\xi_{\tau}}{(x^{\star}-\mathbb{E}(X))(1-\gamma^{-1})}(1+\operatorname{o}(1)) + \frac{\gamma}{\rho(\gamma+\rho)(1-\gamma-\rho)} A(1/\overline{F}(\xi_{\tau}))(1+\operatorname{o}(1)) \right).
\end{multline*}
%
This is immediate from Proposition~\ref{prop:asyexpansion},~\eqref{eqn:asyeqF2}~and~\eqref{eqn:asyeq},
 since $|A|$ is regularly varying with~index~$\rho$.
\end{remark}
\begin{remark}[Comparison with Proposition~3.4 in~\cite{maonghu2015}] Proposition~\ref{prop:asyexpansion} is an extension, tailored to our general semiparametric 
GP setting and extended second-order regular variation assumption $\mathcal{C}_2(\gamma,a,\rho,A)$, of Proposition~3.4 in~\cite{maonghu2015}. The latter result is formulated under a different, nonstandard second-order regular variation condition on $\overline{F}$ when~$X$ belongs to the domain of attraction of a Generalized Extreme Value distribution. It is readily checked by straightforward but tedious calculations that their quantities $c$, $\gamma$, $\rho$ and $A(s)$ respectively correspond to (with the notation of 
Lemma~A.3 in Section~A.1) $C^{1/\gamma}$, $-1/\gamma$, $-\rho/\gamma$ and $-C^{-\rho/\gamma} A(s^{-1/\gamma})/(\gamma(\gamma+\rho))$ of the present paper. In particular, 
when their asymptotic expansion applies, it coincides with ours, but
we lift an unnecessary restriction on the second-order parameter $\rho$ that features in their result.
\end{remark}
\begin{remark}[Expectiles are less extreme than quantiles] An immediate consequence of Equation~\eqref{eqn:asyeq} is that $(x^{\star}-\xi_{\tau})/(x^{\star}-q_{\tau})\to \infty$ as $\tau\uparrow 1$, that is, 
extreme quantiles are closer to the endpoint of a short-tailed distribution than extreme expectiles.
It is therefore unsurprising that the bias due to 
the approximation of tail expectiles by their quantile analogs under the second-order framework,
which is asymptotically proportional to $A((1-\tau)^{-1} (x^{\star} - q_{\tau})^{1/(1-\gamma)})$, converges more slowly to 0 than the corresponding bias term in the heavy-tailed setting, whose order is $A((1-\tau)^{-1})$, see Proposition~1(i) in~\cite{daogirstu2020}. As a second consequence, at least as far as handling bias is concerned, estimating extreme expectiles 
under short-tailed models
using a semiparametric extreme value methodology should be expected to be much harder than 
under heavy-tailed models.
\end{remark}
\begin{remark}[On remainder terms in the asymptotic expansion] The quantity $x^{\star}-\mathbb{E}(X)$, which is a measure of the spread of the distribution tail, appears in the asymptotic equivalent of $(x^{\star} - \xi_{\tau})/(x^{\star} - q_{\tau})$ and in both the remainder terms of the asymptotic expansion for $x^{\star} - \xi_{\tau}$. By contrast, no measures of spread appear in the asymptotic 
connection between extreme expectiles and quantiles
of heavy-tailed distributions, although the expectation $\mathbb{E}(X)$, which can be understood as a location parameter, appears in an error term proportional to $1/q_{\tau}$, as can be seen from 
Proposition~1 in~\cite{daogirstu2020}.
\end{remark}
With Proposition~\ref{prop:asyexpansion} at our disposal, we can now construct and study two classes of extreme expectile estimators. The first one, in Section~\ref{sec:main:LAWS} below, is built upon asymmetric least squares minimization, while the second one, in Section~\ref{sec:main:QB}, 
is directly obtained by plugging in Equation~\eqref{eqn:asyeq} estimators of $\mathbb{E}(X)$ and of the tail quantities $\gamma$, $x^{\star}$ and $q_{\tau}$.

\subsection{Asymmetric least squares estimation}
\label{sec:main:LAWS}

Suppose that the available data 
has been generated from the random variables $X_1,\ldots,X_n$ with common distribution function $F$, and let $\tau_n\uparrow 1$ (as $n\to \infty$) be a high asymmetry level 
at which the target unknown expectile $\xi_{\tau_n}$ is to be estimated.
A first solution is to construct the estimator minimizing the empirical counterpart of problem~\eqref{eqn:expectile}. This produces the Least Asymmetrically Weighted Squares (LAWS) estimator 
\begin{equation}
\label{eq:empirical_LAWS}
\widehat{\xi}_{\tau_n} = \argmin_{\theta\in \mathbb{R}} \frac{1}{n} \sum_{t=1}^n  \eta_{\tau_n}(X_t-\theta) - 
 \eta_{\tau_n} (X_t)  
= \argmin_{\theta\in \mathbb{R}} \sum_{t=1}^n  \eta_{\tau_n}(X_t-\theta).
\end{equation}
A possible way to derive the asymptotic normality of $\widehat{\xi}_{\tau_n}$ would be to find a quadratic approximation to $\sum_{t=1}^n \eta_{\tau_n}(X_t-\theta)$ and then to use stochastic convex optimization theorems such as those of~\cite{hjopol1993}. 
 However, such developments already require lengthy technical arguments in the better-known heavy-tailed model, as can be seen from~\cite{daogirstu2018}. Instead, we
propose to use here an alternative technique leading to a much simpler proof, 
which is
based on the following observation made by~\cite{jon1994}: the $\tau$th expectile of $F$ is 
the $\tau$th quantile of the distribution function $E=1-\overline{E}$,~where
\[
\overline{E}(x) = \frac{\mathbb{E} ( |X-x| \ind{\{ X>x \}} )}{\mathbb{E} ( |X-x| )}. 
\]
This survival function can equivalently be rewritten as 
\[
\overline{E}(x) = \frac{\varphi^{(1)}(x)}{2 \varphi^{(1)}(x) + x - \mathbb{E}(X)}, \mbox{ with } \varphi^{(\kappa)}(x) = \mathbb{E} ( (X-x)^{\kappa} \ind{\{ X>x \}} ).
\]
Since $\widehat{\xi}_{\tau}$ is the $\tau$th expectile of the empirical distribution function 
\[
\widehat{F}_n(x)=\frac{1}{n} \sum_{t=1}^n \ind{\{ X_t>x \}},
\]
it must therefore be the $\tau$th quantile of the distribution function $\widehat{E}_n=1-\widehat{\overline{E}}_n$ defined as 
\[
\widehat{\overline{E}}_n(x) = \frac{\widehat{\varphi}_n^{(1)}(x)}{2 \widehat{\varphi}_n^{(1)}(x) + x - \overline{X}_n}, \mbox{ where } \widehat{\varphi}_n^{(\kappa)}(x) = \frac{1}{n} \sum_{t=1}^n (X_t-x)^{\kappa} \ind{\{ X_t>x \}} ,
\]
with $\overline{X}_n$ being
the sample mean. 
Intuitively, to derive now
the asymptotic behavior of $\widehat{\xi}_{\tau_n}-\xi_{\tau_n}$, it suffices to obtain the asymptotic behavior of $\widehat{\overline{E}}_n(x)/\overline{E}(x)$ at a level $x=x_n$ close to $\xi_{\tau_n}$ in an appropriate sense and to apply a suitable inversion argument.

We do so in a general framework of strictly stationary, weakly dependent random variables. Recall that a strictly stationary sequence $(X_t)_{t\geq 1}$ is said to be $\alpha-$mixing (or strongly mixing) if $\alpha(l)=\sup_{m\geq 1} \alpha_m(l)\to 0$ when $l\to\infty$, where 
\[
\alpha_m(l)=\sup_{\substack{A\in \mathcal{F}_{1,m} \\ B\in \mathcal{F}_{m+l,\infty}}} |\mathbb{P}(A\cap B)-\mathbb{P}(A)\mathbb{P}(B)|
\]
with 
$\mathcal{F}_{1,m}=\sigma(X_1,\ldots,X_m)$ 
and 
$\mathcal{F}_{m+l,\infty}=\sigma(X_{m+l},X_{m+l+1},\ldots)$ 
being
past and future $\sigma-$algebras. The $\alpha-$mixing condition is one of the weakest dependence assumptions in the mixing time series literature: more restrictive conditions include $\beta-$, $\rho-$, $\phi-$ and $\psi-$mixing, see~\cite{bra2005}. We make the following assumption about the mixing rate.
\vskip1ex
\noindent
Condition $\mathcal{M} \ $ There exist sequences of positive integers $(l_n)$ and $(r_n)$, both tending to infinity, such that $l_n/r_n\to 0$, $r_n/n\to 0$ and $n \, \alpha(l_n)/r_n\to 0$, 
as $n\to\infty$.
\vskip1ex
\noindent
The sequences $(l_n)$ and $(r_n)$ are respectively interpreted as ``small-block'' and ``big-block'' sequences, and are used to develop a big-block/small-block argument as a prerequisite to evaluating the asymptotic variance of $\widehat{\xi}_{\tau_n}$. Condition $\mathcal{M}$ has already been used in the literature on the extreme values of time series, 
see {\it e.g.}~\citet[][Equation~(2.1)]{rooleahaa1998}. We also require the 
following tail dependence condition on the joint extreme behavior of $(X_t)_{t\geq 1}$ at different time points.
\vskip1ex
\noindent
Condition $\mathcal{D} \ $ For any integer $t\geq 1$, there is a function $R_t$ on $[0,\infty]^2\setminus \{ (\infty, \infty) \}$ such that 
\[
\forall (x,y)\in [0,\infty]^2\setminus \{ (\infty, \infty) \}, \ \lim_{s\to\infty} s \,  \mathbb{P}( \overline{F}(X_1)\leq x/s, \overline{F}(X_{t+1})\leq y/s ) = R_t(x,y), 
\]
and there exist a constant $K\geq 0$ and a nonnegative summable sequence $(\rho(t))_{t\geq 1}$ such that, for $s$ large enough, 
\[
\forall t\geq 1, \ \forall x,y\in (0,1], \ s \, \mathbb{P}( \overline{F}(X_1)\leq x/s, \overline{F}(X_{t+1})\leq y/s )\leq \rho(t) \sqrt{xy} + \frac{K}{s} xy. 
\]
\vskip1ex
\noindent
The function $R_t$, called the tail copula of $(X_1,X_{t+1})$~\citep[see][]{schsta2006}, finely quantifies the degree of asymptotic dependence between $X_1$ and $X_{t+1}$. Condition~$\mathcal{D}$ ensures that the probability of a joint extreme value of $X_1$ and $X_{t+1}$ is of the same order of magnitude as the probability of an extreme value of $X_1$, meaning that clusters of extreme values across time cannot form too often. A similar anti-clustering assumption is made in~\cite{dre2003}, see conditions (C2) and (C3) therein.

Under these temporal dependence assumptions and using our insight about the link between the LAWS estimator and the empirical estimator of $\overline{E}$, we can prove the following result on the joint asymptotic normality of the LAWS estimator and an empirical quantile having the same order of magnitude.
\begin{theorem}
\label{theo:asynor_LAWS}
Assume that $X$ satisfies condition $\mathcal{C}_2(\gamma,a,\rho,A)$. Let $\tau_n,\alpha_n\uparrow 1$ be such that $n \overline{F}(\xi_{\tau_n})\to\infty$, $\overline{F}(\xi_{\tau_n})/(1-\alpha_n)\to 1$ and $\sqrt{n \overline{F}(\xi_{\tau_n})} A(1/\overline{F}(\xi_{\tau_n}))=\operatorname{O}(1)$.
\begin{enumerate}[label=(\roman*)]
\item Suppose that $(X_t)_{t\geq 1}$ is a strictly stationary sequence of copies of $X$, whose distribution function $F$ is continuous, satisfying conditions $\mathcal{M}$ and $\mathcal{D}$. Assume that $r_n \overline{F}(\xi_{\tau_n})\to 0$, and that there is $\delta>0$ such that
\[
\mathbb{E}(|\min(X,0)|^{2+\delta})<\infty, \ \sum_{l\geq 1} l^{2/\delta} \alpha(l)<\infty \mbox{ and } r_n\left(\frac{r_n}{\sqrt{n\overline{F}(\xi_{\tau_n})}} \right)^{\delta}\to 0. 
\]
Then 
\[
\frac{\sqrt{n \overline{F}(\xi_{\tau_n})}}{a(1/\overline{F}(\xi_{\tau_n}))}  (\widehat{\xi}_{\tau_n}-\xi_{\tau_n}, \widehat{q}_{\alpha_n}-q_{\alpha_n}) \stackrel{d}{\longrightarrow} \mathcal{N}( 0, \bm{V}(\gamma)+2\,\bm{C}(\gamma,R) )
\]
where the $2\times 2$ symmetric matrices $\bm{V}(\gamma)$ and $\bm{C}(\gamma,R)$ are defined elementwise as $V_{11}(\gamma)=2/[(1-\gamma)(1-2\gamma)]$, $V_{12}(\gamma)=1/(1-\gamma)$ and $V_{22}(\gamma)=1$, 
\begin{align*}
C_{11}(\gamma,R) &= \frac{1}{\gamma^2} \iint_{(0,1]^2} \sum_{t=1}^{\infty} R_t(x^{-1/\gamma},y^{-1/\gamma}) \, \mathrm{d}x \, \mathrm{d}y \\
C_{12}(\gamma,R) &= -\frac{1}{2\gamma} \int_0^1 \sum_{t=1}^{\infty} [R_t(x^{-1/\gamma},1) + R_t(1,x^{-1/\gamma})] \, \mathrm{d}x \\ 
\mbox{and } C_{22}(\gamma,R) &= \sum_{t=1}^{\infty} R_t(1,1).
\end{align*}
\item If the $X_i$ are i.i.d.~copies of $X$ and $\mathbb{E}(|\min(X,0)|^2)<\infty$, then the above asymptotic normality result holds with $R_t\equiv 0$ for any $t\geq 1$, that is, 
\[
\frac{\sqrt{n \overline{F}(\xi_{\tau_n})}}{a(1/\overline{F}(\xi_{\tau_n}))}  (\widehat{\xi}_{\tau_n}-\xi_{\tau_n}, \widehat{q}_{\alpha_n}-q_{\alpha_n}) \stackrel{d}{\longrightarrow} \mathcal{N}( 0, \bm{V}(\gamma) ).
\]
\end{enumerate}
If $X$ is bounded, then assumption $\sum_{l\geq 1} l^{2/\delta} \alpha(l)<\infty$ in (i) can be weakened to $\sum_{l\geq 1} \alpha(l)<\infty$, and no integrability assumption on $X$ is necessary.
\end{theorem}
In Theorem~\ref{theo:asynor_LAWS}, condition $n \overline{F}(\xi_{\tau_n})\to\infty$ requires that 
$\tau_n$ be intermediate, {\it i.e.}~not too large. 
Assumption
$\sqrt{n \overline{F}(\xi_{\tau_n})} A(1/\overline{F}(\xi_{\tau_n}))=\operatorname{O}(1)$ is a bias condition which corresponds exactly
to the usual bias condition $\sqrt{n (1-\tau_n)} A((1-\tau_n)^{-1})=\operatorname{O}(1)$ in extreme quantile estimation when replacing $\xi_{\tau_n}$ 
by $q_{\tau_n}$,
see Theorem~2.4.1 on p.50 of~\cite{haafer2006}. The conditions on $r_n$ in Theorem~\ref{theo:asynor_LAWS}(i) are similar to those of Theorem~3.1 in~\cite{davpadstu2022} 
under heavy-tailed models,
taking into account that $\overline{F}(\xi_{\tau_n})$ is asymptotically proportional to $1-\tau_n$ in the latter setting. The integrability assumption on $X$ and the 
condition
on the mixing rate $\alpha(l)$ ensure that a central limit theorem applies to $\overline{X}_n$, as part of the proof of the asymptotic normality of $\widehat{\overline{E}}_n(x)/\overline{E}(x)$ at high levels $x=x_n$ close to $\xi_{\tau_n}$.
\begin{remark}[On bias conditions] An inspection of the proof of Theorem~\ref{theo:asynor_LAWS}(i) reveals that the bias condition $\sqrt{n \overline{F}(\xi_{\tau_n})} A(1/\overline{F}(\xi_{\tau_n}))=\operatorname{O}(1)$ is only 
needed for
the asymptotic normality of $\widehat{q}_{\alpha_n}-q_{\alpha_n}$, and is 
unnecessary for the validity of the asymptotic normality of $\widehat{\xi}_{\tau_n}-\xi_{\tau_n}$~alone.
\end{remark}
\begin{remark}[Comparison with the i.i.d.~case] It is natural to compare the asymptotic normality of the LAWS estimator for short-tailed data 
with
the corresponding result one obtains for heavy-tailed data. 
We restrict the comparison to
the i.i.d.~setting for the sake of simplicity. If $X$ has a heavy right tail, that is, $U(sz)/U(s)\to z^{\gamma}$ as $s\to\infty$, for any $z>0$, where $0<\gamma<1/2$, and under the assumptions that $\mathbb{E}(|\min(X,0)|^{2+\delta})<\infty$ for some $\delta>0$, $\tau_n\uparrow 1$ and $n(1-\tau_n)\to\infty$, one has, by Theorem~2 in~\cite{daogirstu2018},
\[
\sqrt{n(1-\tau_n)}\left( \frac{\widehat{\xi}_{\tau_n}}{\xi_{\tau_n}} - 1 \right) \stackrel{d}{\longrightarrow} \mathcal{N}\left( 0, \frac{2\gamma^3}{1-2\gamma} \right).
\]
In this same setting, $(U(sz)-U(s))/a(s) \to (z^{\gamma}-1)/\gamma$ as $s\to\infty$, with $a(s)=\gamma U(s)$, and $\overline{F}(\xi_{\tau_n})/(1-\tau_n)\to \gamma^{-1}-1=(1-\gamma)/\gamma$~\citep[this was first shown by][Theorem~11]{belklamulgia2014}. Therefore, when $X$ has a heavy right tail, 
\[
\frac{\sqrt{n \overline{F}(\xi_{\tau_n})}}{a(1/\overline{F}(\xi_{\tau_n}))}  (\widehat{\xi}_{\tau_n}-\xi_{\tau_n}) \approx \frac{\sqrt{1-\gamma}}{\gamma^{3/2}} \times \sqrt{n (1-\tau_n)} \left( \frac{\widehat{\xi}_{\tau_n}}{\xi_{\tau_n}} - 1 \right) \stackrel{d}{\longrightarrow} \mathcal{N}\left( 0, \frac{2(1-\gamma)}{1-2\gamma} \right).
\]
It follows that the rates of convergence of the LAWS estimator look similar in 
both
heavy and bounded tail settings, but there is a phase transition in terms of asymptotic variance: 
the term $1-\gamma$ appears in its numerator for heavy tails, while it appears in the denominator for short tails, as established in Theorem~\ref{theo:asynor_LAWS}. Interestingly, the two asymptotic variances in the heavy and short-tailed settings converge to 2, and therefore exactly match together in the light-tailed middle scenario, when $\gamma\to 0$.
\end{remark}
\begin{remark}[Comparison with the asymptotic normality of intermediate quantiles] It is instructive to compare Theorem~\ref{theo:asynor_LAWS} with 
the asymptotic normality result for the direct intermediate quantile estimator at level $\tau_n$, 
namely, the empirical counterpart $\widehat{q}_{\tau_n}=X_{\lceil n\tau_n \rceil,n}$.
According to Theorem~2.4.1 on p.50 of~\cite{haafer2006}, when the $X_i$ are i.i.d.,
\[
\frac{\sqrt{n \overline{F}(q_{\tau_n})}}{a(1/\overline{F}(q_{\tau_n}))}   (\widehat{q}_{\tau_n}-q_{\tau_n}) \stackrel{d}{\longrightarrow} \mathcal{N}(0,1).
\]
Observe that, by a combination of Lemma~1.2.9 in~\cite{haafer2006} and 
Lemma~A.1 in Section~A.1, $a(1/\overline{F}(x))/(x^{\star}-x)\to -\gamma$ as $x\uparrow x^{\star}$, and therefore
\[
\left. \frac{\sqrt{n \overline{F}(\xi_{\tau_n})}}{a(1/\overline{F}(\xi_{\tau_n}))} \right/ \frac{\sqrt{n \overline{F}(q_{\tau_n})}}{a(1/\overline{F}(q_{\tau_n}))} = \sqrt{\frac{\overline{F}(\xi_{\tau_n})}{\overline{F}(q_{\tau_n})}} \times \frac{x^{\star}-q_{\tau_n}}{x^{\star}-\xi_{\tau_n}} (1+\operatorname{o}(1)). 
\]
By~\eqref{eqn:asyeqF2} and~\eqref{eqn:asyeq} this ratio is asymptotically proportional to $(x^{\star}-q_{\tau_n})^{-(\gamma+1/2)/(1-\gamma)}$ under the mild further condition $\rho<0$. In other words, the intermediate LAWS estimator $\widehat{\xi}_{\tau_n}$ converges faster than $\widehat{q}_{\tau_n}$ when $\gamma>-1/2$, has the same rate of convergence if $\gamma=-1/2$, and converges at a slower rate if $\gamma<-1/2$.
\end{remark}
As a corollary of Theorem~\ref{theo:asynor_LAWS}, we obtain the asymptotic normality of the empirical estimator $\widehat{\overline{F}}_n(\widehat{\xi}_{\tau_n})$ of $\overline{F}(\xi_{\tau_n})$, on which the rate of convergence of $\widehat{\xi}_{\tau_n}$ crucially depends.
\begin{corollary}
\label{coro:asynor_surv}
Work under the conditions of Theorem~\ref{theo:asynor_LAWS}. Then
\begin{multline*}
\sqrt{n \overline{F}(\xi_{\tau_n})} \left( \frac{\widehat{\overline{F}}_n(\widehat{\xi}_{\tau_n})}{\overline{F}(\xi_{\tau_n})}-1 \right) \stackrel{d}{\longrightarrow} \mathcal{N}\left( 0, \frac{2\gamma^2+\gamma+1}{(1-\gamma)(1-2\gamma)} \right. \\
\left. + 2\iint_{(0,1]^2}\sum_{t=1}^{\infty} \left( \frac{1}{\gamma^2} R_t(x^{-1/\gamma},y^{-1/\gamma}) + \frac{1}{\gamma} [R_t(x^{-1/\gamma},1) + R_t(1,x^{-1/\gamma})] + R_t(1,1) \right) \mathrm{d}x \, \mathrm{d}y \right).
\end{multline*}
\end{corollary}
\begin{remark}[On rates of convergence] The rate of convergence of $\widehat{\overline{F}}_n(\widehat{\xi}_{\tau_n})$ is rather natural: for a sequence $(u_n)$ tending to $x^{\star}$ such that $n\overline{F}(u_n)\to \infty$, 
Lemma~A.5 states that 
\[
\sqrt{n \overline{F}(u_n)} \left( \frac{\widehat{\overline{F}}_n(u_n)}{\overline{F}(u_n)} - 1 \right) \stackrel{d}{\longrightarrow} \mathcal{N}\left( 0, 1+2\sum_{t=1}^{\infty} R_t(1,1) \right).
\]
It is worth noticing that the asymptotic variance of $\widehat{\overline{F}}_n(\widehat{\xi}_{\tau_n})$ does not coincide with the variance that would be obtained if $\xi_{\tau_n}$ were known, namely, if $\widehat{\overline{F}}_n(\xi_{\tau_n})$ were considered instead. This is due to the asymptotic dependence existing between $\widehat{\xi}_{\tau_n}$ and high order statistics (and therefore between $\widehat{\xi}_{\tau_n}$ and $\widehat{\overline{F}}_n$), see Theorem~\ref{theo:asynor_LAWS} and the proof of Corollary~\ref{coro:asynor_surv}.
\end{remark}
We now have the tools necessary to construct an extreme value estimator of a properly extreme expectile $\xi_{1-p_n}$, where $p_n\downarrow 0$ at any possible rate as $n\to\infty$, by extrapolating the intermediate LAWS estimator $\widehat{\xi}_{\tau_n}$ to the right place at the far tail. Using~\eqref{eqn:approx} with $s=1/\overline{F}(\xi_{\tau_n})$ and $z=\overline{F}(\xi_{\tau_n})/\overline{F}(\xi_{1-p_n})$ motivates the approximation
\[
\xi_{1-p_n} \approx \xi_{\tau_n} + a(1/\overline{F}(\xi_{\tau_n})) \frac{(\overline{F}(\xi_{\tau_n})/\overline{F}(\xi_{1-p_n}))^{\gamma}-1}{\gamma}.
\]
By Theorem~\ref{theo:asynor_LAWS}, $\xi_{\tau_n}$ is estimated by the LAWS estimator $\widehat{\xi}_{\tau_n}$ at rate $a(1/\overline{F}(\xi_{\tau_n}))/\sqrt{n \overline{F}(\xi_{\tau_n})}$. The scale parameter $a(1/\overline{F}(\xi_{\tau_n}))$ and shape parameter $\gamma$ can be estimated by a variety of techniques, such as the 
GP-pseudo-ML 
estimator and 
Moment-type
estimator, see Sections~3.4 and~3.5 in~\cite{haafer2006}. Typical such estimators of the scale function $a(1/\overline{F}(u_n))$, 
when
$u_n\uparrow x^{\star}$ 
such that $n \overline{F}(u_n)\to\infty$, converge on the relative scale at the rate $1/\sqrt{n \overline{F}(u_n)}$; see Sections~3.4 and 4.2 in~\cite{haafer2006} in the i.i.d.~case, and Section~6 in~\cite{dre2003} in the dependent data setup. Since, by Corollary~\ref{coro:asynor_surv}, the (unknown) quantity $\overline{F}(u_n)=\overline{F}(\xi_{\tau_n})$ can be consistently estimated at the rate $1/\sqrt{n \overline{F}(\xi_{\tau_n})}$, we therefore expect to be able to estimate $a(1/\overline{F}(\xi_{\tau_n}))$ at this rate on the relative scale. Finally, given an intermediate level $\tau_n$, it is customary to estimate the extreme value index $\gamma$ at the rate $1/\sqrt{n(1-\tau_n)}$ when the top $k=\lfloor n(1-\tau_n) \rfloor$ values in the data are used, see Sections~3.4, 3.5 and 3.6 in~\cite{haafer2006} in the i.i.d.~case, and again Section~6 in~\cite{dre2003} when the data points are serially dependent. It remains to find a way to estimate $\overline{F}(\xi_{\tau_n})/\overline{F}(\xi_{1-p_n})$, which depends on the target quantity $\xi_{1-p_n}$ itself. 
Combining~\eqref{eqn:asyeqF2} and~\eqref{eqn:asyeq} with the fact that the function $s\mapsto x^{\star}-U(s)$ is regularly varying with index $\gamma$~\citep[][Corollary~1.2.10 on~p.23]{haafer2006} suggests that 
\begin{align}
\nonumber
    \frac{\overline{F}(\xi_{\tau_n})}{\overline{F}(\xi_{1-p_n})} \approx \frac{1-\tau_n}{p_n} \times \frac{x^{\star}-\xi_{1-p_n}}{x^{\star}-\xi_{\tau_n}} &\approx \frac{1-\tau_n}{p_n} \left( \frac{x^{\star}-q_{1-p_n}}{x^{\star}-q_{\tau_n}} \right)^{1/(1-\gamma)} \\
\label{eqn:ratioFbar}
    &\approx \frac{1-\tau_n}{p_n} \left( \frac{1-\tau_n}{p_n} \right)^{\gamma/(1-\gamma)} = \left( \frac{1-\tau_n}{p_n} \right)^{1/(1-\gamma)} 
\end{align}
%
which in turn leads to the expectile-specific 
approximation
\[
\xi_{1-p_n} \approx \xi_{\tau_n} + a(1/\overline{F}(\xi_{\tau_n})) \frac{((1-\tau_n)/p_n)^{\gamma/(1-\gamma)}-1}{\gamma}.
\]
Consequently, like extreme quantiles, extreme expectiles can be extrapolated from their values at lower levels. More specifically, 
given estimators $\widehat{\sigma}_n$ and $\widehat{\gamma}_n$ of $a(1/\overline{F}(\xi_{\tau_n}))$ and $\gamma$, respectively, one can then construct the 
$\xi_{1-p_n}$ estimator
\begin{equation}\label{eq:extrapo_LAWS}
\widehat{\xi}_{1-p_n}^{\star} = \widehat{\xi}_{\tau_n} + \widehat{\sigma}_n \frac{((1-\tau_n)/p_n)^{\widehat{\gamma}_n/(1-\widehat{\gamma}_n)}-1}{\widehat{\gamma}_n}.
\end{equation}
Since $(1-\tau_n)/\overline{F}(\xi_{\tau_n})\to 0$, the parameter $\gamma$ is estimated at a slower rate than the other quantities, so we expect the asymptotic behavior of $\widehat{\gamma}_n$ to 
govern that of $\widehat{\xi}_{1-p_n}^{\star}$.
The 
last
theorem of this section makes this intuition rigorous. Its proof crucially 
relies on
Theorem~\ref{theo:asynor_LAWS} 
and on Proposition~\ref{prop:asyexpansion} in order to quantify the bias in 
the approximation~\eqref{eqn:ratioFbar}. 
\begin{theorem}
\label{theo:asyext}
Work under the conditions of Theorem~\ref{theo:asynor_LAWS}. If moreover $\rho<0$, $n(1-\tau_n)\to\infty$, $(1-\tau_n)/p_n\to\infty$, $\sqrt{n(1-\tau_n)}/\log((1-\tau_n)/p_n)\to\infty$, $\sqrt{n(1-\tau_n)} (x^{\star}-q_{\tau_n})^{1/(1-\gamma)}=\operatorname{O}(1)$, $\sqrt{n(1-\tau_n)} A((1-\tau_n)^{-1})=\operatorname{O}(1)$, 
$\widehat{\sigma}_n$ and $\widehat{\gamma}_n$ are such that 
\[
\sqrt{n \overline{F}(\xi_{\tau_n})} \left( \frac{\widehat{\sigma}_n}{a(1/\overline{F}(\xi_{\tau_n}))} - 1 \right) = \operatorname{O}_{\mathbb{P}}(1) \ \mbox{ and } \ \sqrt{n(1-\tau_n)} (\widehat{\gamma}_n-\gamma)\stackrel{d}{\longrightarrow} \Gamma,
\]
where $\Gamma$ is a nondegenerate limit, then
\[
\frac{\sqrt{n(1-\tau_n)}}{a(1/\overline{F}(\xi_{\tau_n}))} (\widehat{\xi}_{1-p_n}^{\star} - \xi_{1-p_n}) \stackrel{d}{\longrightarrow} \frac{\Gamma}{\gamma^2}.
\]
\end{theorem}
\subsection{Quantile-based estimation}
\label{sec:main:QB}

We use here Proposition~\ref{prop:asyexpansion} to 
present an alternative estimator of extreme expectiles, purely based on quantiles, and to develop its asymptotic theory.
Similarly to the setup of extreme quantile estimation in
Section~4.3 of~\cite{haafer2006}, assume that $k=k_n\to\infty$ is a sequence of positive integers such that $k/n\to 0$ and that respective estimators $\widehat{\gamma}_n$, $\widehat{a}(n/k)$ and $X_{n-k,n}$ of $\gamma$, $a(n/k)$ and $U(n/k)$ 
are given such that 
\begin{equation}
\label{eq:basic_conditions}
\sqrt{k}\left( \widehat{\gamma}_n-\gamma, \frac{\widehat{a}(n/k)}{a(n/k)}-1, \frac{X_{n-k,n}-U(n/k)}{a(n/k)} \right) \stackrel{d}{\longrightarrow} (\Gamma, \Lambda, B)
\end{equation}
where $(\Gamma, \Lambda, B)$ is a nontrivial trivariate weak limit. This assumption is satisfied by moment and 
GP-pseudo-ML
estimators of the 
shape
and scale parameters, among others, see an overview in Section~4.2 of~\cite{haafer2006} in the case 
where
the $X_i$ are independent random variables. It is also satisfied when $(X_t)_{t\geq 1}$ is a strictly stationary but serially dependent sequence: this is for example the case when the data points are $\beta-$mixing and satisfy an anti-clustering condition similar to the tail dependence assumption $\mathcal{D}$, as a consequence of the powerful results of~\cite{dre2003}. 

Let $p_n\downarrow 0$ 
with $k/(n p_n)\to\infty$, so that the level $1-p_n$
is much more extreme than $1-k/n$. 
Following Section~4.3 in~\cite{haafer2006}, 
the extreme quantile $q_{1-p_n}$ and the right endpoint $x^{\star}$ can be estimated by 
\begin{equation}
\label{eq:quantile_endpoint}
\widehat{q}_{1-p_n}^{\star} = X_{n-k,n} + \widehat{a}(n/k) \frac{(k/(np_n))^{\widehat{\gamma}_n}-1}{\widehat{\gamma}_n} \ \mbox{ and } \ \widehat{x}^{\star} = X_{n-k,n} - \frac{\widehat{a}(n/k)}{\widehat{\gamma}_n}.
\end{equation}
According to Proposition~\ref{prop:asyexpansion}, an estimator of $\xi_{1-p_n}$ 
follows then as
\begin{equation}
\label{eq:extrapo_QB}
\widetilde{\xi}_{1-p_n}^{\star} = \widehat{x}^{\star} - [(\widehat{x}^{\star}-\overline{X}_n)(1-\widehat{\gamma}_n^{-1})]^{-\widehat{\gamma}_n/(1-\widehat{\gamma}_n)} (\widehat{x}^{\star} - \widehat{q}_{1-p_n}^{\star})^{1/(1-\widehat{\gamma}_n)}.
\end{equation}
The next result provides its asymptotic properties, where 
two sequences $(u_n)$ and $(v_n)$ are said to be asymptotically proportional if $(u_n/v_n)$ 
tends to a finite positive limit as $n\to\infty$. 
\begin{theorem}
\label{theo:asynor_QB}
Suppose that $\mathbb{E}|\min(X,0)|<\infty$ and condition $\mathcal{C}_2(\gamma,a,\rho,A)$ holds with $\rho<0$, and let $x^{\star}$ be the finite right endpoint of $F$. 
Assume that condition \eqref{eq:basic_conditions} holds true and that $\sqrt{k}(\overline{X}_n-\mathbb{E}(X))\stackrel{\mathbb{P}}{\longrightarrow} 0$
%
with $k=k_n$ 
being
asymptotically proportional to $n^{\alpha}$, for some $\alpha\in (0,1)$.
Let $p_n$ be asymptotically proportional to $n^{-\beta}$ where $\beta>0$ is such that $\alpha+\beta-1>0$. If moreover $\sqrt{k} A(n/k)\to \lambda\in \mathbb{R}$, then we have, up to changing probability spaces and with appropriate versions of the estimators involved,
\begin{align*}
    \widetilde{\xi}_{1-p_n}^{\star} - \xi_{1-p_n} &= \frac{a(\frac{n}{k})}{\sqrt{k}} \frac{1}{\gamma^2} \left( \Gamma + \gamma^2 B - \gamma \Lambda - \frac{\lambda\gamma}{\gamma+\rho} + \operatorname{o}_{\mathbb{P}}(1) \right) \\
    &- \left[a\left(\frac{n}{k}\right)\left(\frac{k}{n p_n}\right)^{\gamma}\right]^{\frac{1}{1-\gamma}}  \frac{(1-\gamma)^{-\frac{1}{1-\gamma}}}{\gamma (x^{\star}-\mathbb{E}(X))^{\frac{\gamma}{1-\gamma}} } 
    \left( \frac{\log(\frac{n p_n^{1/(1-\gamma)}}{k})}{\sqrt{k}} \Gamma +\operatorname{o}_{\mathbb{P}}\left( \frac{\log n}{\sqrt{k}} \right) 
    \right) \\
    &+ \operatorname{O}(n^{\beta\gamma/(1-\gamma)} ( n^{\beta\gamma/(1-\gamma)} + |A( n^{\beta/(1-\gamma)} )| )).
\end{align*} 
\end{theorem}
\begin{remark}[Mixing and the central limit theorem] Condition $\sqrt{k}(\overline{X}_n-\mathbb{E}(X))\stackrel{\mathbb{P}}{\longrightarrow} 0$ 
is satisfied in practice
if $\sqrt{n}(\overline{X}_n-\mathbb{E}(X))=\operatorname{O}_{\mathbb{P}}(1)$, which is in particular true when a central limit theorem applies. As already highlighted below Theorem~\ref{theo:asynor_LAWS}, this will be the case if $\mathbb{E}(|\min(X,0)|^2)<\infty$ when the $X_i$ are independent, or if there is $\delta>0$ such that $\mathbb{E}(|\min(X,0)|^{2+\delta})<\infty$ and $\sum_{l\geq 1} l^{2/\delta} \alpha(l)<\infty$ when $(X_t)_{t\geq 1}$ is $\alpha-$mixing. In particular, when the data mixes geometrically fast, then $\sqrt{n}(\overline{X}_n-\mathbb{E}(X))=\operatorname{O}_{\mathbb{P}}(1)$ as soon as $X$ has a finite moment of order $2+\delta$, for some $\delta>0$.
\end{remark}
%
\begin{remark}[Our assumptions on $k$] The assumption that $k=k_n$ is asymptotically equivalent to a positive and finite multiple of $n^{\alpha}$, is only very slightly stronger than the usual pair of extreme value conditions $k\to\infty$ and $k/n\to 0$ made throughout Section~4 in~\cite{haafer2006}. 
The only difference is that 
our assumption does not allow to take $k$ growing to infinity logarithmically fast; 
such sequences produce, however, 
very small values of $k$ in practice and would therefore yield estimators having very large variances. We also note that in standard settings such as those of~\citet[][Table~2.2 on~p.68]{beigoesegteu2004}, $A(s)$ is asymptotically proportional to $s^{\rho}$, in which case the optimal choices of $k$ satisfying the usual bias-variance tradeoff for extreme value index estimation would fulfill $\sqrt{k} A(n/k)\to \lambda\in \mathbb{R}\setminus \{ 0 \}$, that is, 
$k=\operatorname{O}(n^{-2\rho/(1-2\rho)})$.
In other words, it is reasonable to expect that optimal choices of $k$ in practice have to be asymptotically equivalent to a positive and finite multiple of a fractional power of $n$.
\end{remark}
It follows from Theorem~\ref{theo:asynor_QB} that the asymptotic behavior of the extreme expectile estimator $\widetilde{\xi}_{1-p_n}^{\star}$ is more complex than that of the extreme quantile estimator $\widehat{q}_{1-p_n}^{\star}$: while, from Theorem~4.3.1 on p.134 and Theorem~4.5.1 on~p.146 of~\cite{haafer2006}, 
$\widehat{q}_{1-p_n}^{\star}-q_{1-p_n}$ converges to the same distribution $\frac{1}{\gamma^2} \left( \Gamma + \gamma^2 B - \gamma \Lambda - \lambda \frac{\gamma}{\gamma+\rho} \right)$ as $\widehat{x}^{\star} - x^{\star}$
at the rate $a(n/k)/\sqrt{k}$ 
for
$\gamma<0$, the asymptotic distribution of $\widetilde{\xi}_{1-p_n}^{\star}-\xi_{1-p_n}$ may be a nonstandard mixture of 
the two limiting distributions of $\widehat{x}^{\star} - x^{\star}$ and $\widehat{\gamma}_n-\gamma$. In particular,
Corollary~\ref{coro:asynor} shows that 
when, for example,
$\beta=1$ and $\alpha$ (and hence $k$) is chosen small enough, it is in fact 
the asymptotic distribution $\Gamma$ of $\widehat{\gamma}_n-\gamma$
that dominates in $\widetilde{\xi}_{1-p_n}^{\star}-\xi_{1-p_n}$, 
while
Corollary~\ref{coro:asynor_2} 
examines what can otherwise be said.
\begin{corollary}
\label{coro:asynor}
Work under Theorem~\ref{theo:asynor_QB}. If 
$\alpha<1-\beta/(1-\gamma)$ and $\alpha<2\beta\min(-\gamma,-\rho)/(1-\gamma)$, 
%
\begin{multline*}
    \frac{\sqrt{k}}{\log(n p_n^{1/(1-\gamma)}/k)} \frac{\widetilde{\xi}_{1-p_n}^{\star} - \xi_{1-p_n}}{[a(\frac{n}{k}) (\frac{k}{n p_n})^{\gamma} ]^{1/(1-\gamma)}} 
    \stackrel{d}{\longrightarrow} -\gamma^{-1} (1-\gamma)^{-1/(1-\gamma)} (x^{\star}-\mathbb{E}(X))^{-\gamma/(1-\gamma)} \Gamma.
\end{multline*} 
\end{corollary}
\begin{remark}[Link between bias terms and short-tailedness] The closer $\gamma$ is to 0, the stronger the constraint on $\alpha$ through the condition $\alpha<2\beta\min(-\gamma,-\rho)/(1-\gamma)$. This is analogous to what 
happens
in extreme expectile estimation for heavy-tailed distributions, where the condition $\sqrt{k}/q_{1-k/n}=\operatorname{O}(1)$~\citep[see {\it e.g.}][Theorem~5]{daogirstu2020} becomes a strong restriction as the tail gets less heavy, {\it i.e.}~when $\gamma$ approaches 0.
\end{remark}
\begin{remark}[Comparison between the LAWS and quantile-based estimators] 
\label{rmk:rates}
One may compare the rates of convergence of $\widehat{\xi}_{1-p_n}^{\star}$ and $\widetilde{\xi}_{1-p_n}^{\star}$ by setting $\tau_n=1-k/n$. Using the convergence $a(s)/(x^{\star} - U(s))\to -\gamma$ as $s\to \infty$ and Equation~\eqref{eqn:asyeq}, 
one finds under the assumptions of Corollary~\ref{coro:asynor} that
\begin{multline*}
    \frac{\sqrt{k}}{a(1/\overline{F}(\xi_{1-k/n}))} \left/ \frac{\sqrt{k}}{\log(n p_n^{1/(1-\gamma)}/k) [a(n/k) (k/(n p_n))^{\gamma} ]^{1/(1-\gamma)}} \right. \\
    \propto \log(n) \frac{[a(n/k)]^{1/(1-\gamma)}}{a(1/\overline{F}(\xi_{1-k/n}))} (k/(n p_n))^{\gamma/(1-\gamma)} \propto n^{(\alpha+\beta-1)\gamma/(1-\gamma)} \log(n) \to 0. 
\end{multline*}
This means that $\widetilde{\xi}_{1-p_n}^{\star}$ converges to $\xi_{1-p_n}$ faster than $\widehat{\xi}_{1-p_n}^{\star}$ when $k$ 
(or $1-\tau_n$) 
is chosen small. We shall illustrate 
this finding below
in our simulation study.
\end{remark}
%
Condition $\alpha<1-\beta/(1-\gamma)$ may not hold in a given example, especially when $\beta$ is large enough, or equivalently, $p_n$ is small enough. Yet, interestingly this condition can always be satisfied for sufficiently small $\alpha$
in the standard setting $\beta=1$ of extreme value analysis. If 
it is not satisfied, then $\widetilde{\xi}_{1-p_n}^{\star}$ tends to inherit the asymptotic behavior of $\widehat{x}^{\star}$, rather than $\widehat{\gamma}$, as established in the following result.
\begin{corollary}
\label{coro:asynor_2}
Under the assumptions of Theorem~\ref{theo:asynor_QB}, if moreover $\alpha>1-\beta/(1-\gamma)$, then 
%
%
\begin{multline*}
    \widetilde{\xi}_{1-p_n}^{\star} - \xi_{1-p_n} = \frac{a(\frac{n}{k})}{\sqrt{k}\gamma^2}  \left( \Gamma + \gamma^2 B - \gamma \Lambda -  \frac{\lambda\gamma}{\gamma+\rho} + \operatorname{o}_{\mathbb{P}}(1) \right)
  +  \operatorname{O}(n^{\frac{\beta\gamma}{1-\gamma}} ( n^{\frac{\beta\gamma}{1-\gamma}} + |A( n^{\frac{\beta}{1-\gamma}} )| )).
\end{multline*} 
\end{corollary}
The condition $\alpha>1-\beta/(1-\gamma)$ itself is not sufficient to ensure the convergence of $\widetilde{\xi}_{1-p_n}^{\star}$; in practice, the bias term may dominate the asymptotics depending on the choice of $k$. This is most easily seen when $A(s)$ is asymptotically proportional to $s^{\rho}$ and $\beta=1$, corresponding to the standard extreme value situation where $p_n\approx c/n$. In this case, 
(i) one automatically has $\alpha+\beta-1=\alpha>0$, 
(ii) condition $\sqrt{k} A(n/k)\to \lambda\in \mathbb{R}$ essentially amounts to $\alpha\leq -2\rho/(1-2\rho)$, 
and (iii) condition $\alpha>1-\beta/(1-\gamma)$ becomes $\alpha>-\gamma/(1-\gamma)$. 
For the bias term in Corollary~\ref{coro:asynor_2} to be negligible, one requires 
\[
\frac{\sqrt{k}}{a(n/k)} \times n^{\gamma/(1-\gamma)} ( n^{\gamma/(1-\gamma)} + |A( n^{1/(1-\gamma)} )| ) \to 0.
\]
Since $a(n/k)$ is asymptotically proportional to $(n/k)^{\gamma}$ by 
Lemma~A.3(i), this is equivalent to assuming 
\[
\alpha\left( \frac{1}{2} + \gamma \right) + \frac{\gamma^2 - \min(-\gamma,-\rho)}{1-\gamma} < 0. 
\]
When $\gamma>-1/2$, as is often the case in applications, and $0<-\rho<-\gamma$, representing situations where the bias due to the second-order framework is high, this condition becomes 
\[
\alpha< -\frac{2(\gamma^2 + \rho)}{(1-\gamma)(1+2\gamma)}.
\]
Depending on the value of $\rho$, this final condition may not be compatible with $\alpha>-\gamma/(1-\gamma)$: 
in fact, if $\rho$ is close enough to 0, it may even be impossible to satisfy whatever the value of $\alpha$ (since the right-hand side of the above displayed inequality tends to a negative constant as $\rho\to 0$, when $\gamma>-1/2$). In this case, with the choice $p_n=c/n$, the asymptotic behavior of $\widehat{x}^{\star}-x^{\star}$ can never dominate in $\widetilde{\xi}_{1-p_n}^{\star} - \xi_{1-p_n}$.
\subsection{Selection of the expectile asymmetry level}
\label{sec:main:select}

In practical situations it is crucial to make an informed decision as to what the asymmetry level of the target expectile should be. In financial applications, where the dual interpretation of expectiles in terms of the gain-loss ratio is available~\citep{beldib2017}, it is sensible to set 
the expectile level so as
to achieve a certain value of the gain-loss ratio. Otherwise, it has been proposed in the literature to select $\tau$ such that $\xi_{\tau}$ coincides with another pre-specified
intuitive 
risk measure:~\cite{beldib2017} suggest to choose the expectile level $\tau$ so that $\xi_{\tau}$ is 
identical
to the Value-at-Risk (or quantile) $q_{\alpha}$, where $\alpha$ is a high 
tail probability level specified by the statistician or the practitioner.

The proposal of~\cite{beldib2017} is valid only when the underlying loss distribution is Gaussian. \cite{daogirstu2018} later extended this idea to the heavy-tailed setup.
We examine here the short-tailed situation, hitherto unexplored from this perspective. Fix 
a large quantile level $\alpha=\alpha_n=1-p_n$.
Setting $\tau=\tau_n=1-\pi_n$ to be the expectile level such that $\xi_{\tau}=q_{\alpha}$, Equation~\eqref{eqn:asyeqF2} leads to
\[
\frac{(x^{\star} - q_{1-p_n}) p_n}{\pi_n} = \frac{(x^{\star} - \xi_{\tau}) \overline{F}(\xi_{\tau})}{1-\tau} \approx (x^{\star}-\mathbb{E}(X))(1-\gamma^{-1}).
\]
In other words, 
\[
\pi_n \approx \frac{x^{\star} - q_{1-p_n}}{(x^{\star}-\mathbb{E}(X))(1-\gamma^{-1})} p_n.
\]
This approximation suggests to estimate the quantity $\pi_n$ by 
\[
\widehat{\pi}_n\equiv\widehat{\pi}_n(p_n)=\frac{\widehat{x}^{\star} - \widehat{q}_{1-p_n}^{\star}}{(\widehat{x}^{\star}-\overline{X}_n)(1-\widehat{\gamma}_n^{-1})} p_n
\]
with the notation of~\eqref{eq:quantile_endpoint}. Our final main result shows that this estimator is asymptotically normal in the framework of Section~\ref{sec:main:QB}. 
\begin{proposition}
\label{prop:estim_pin}
Work under Theorem~\ref{theo:asynor_QB}. If moreover $\alpha<\min(-2\beta\gamma,-2\rho/(1-2\rho))$, then 
\[
\frac{\sqrt{k}}{\log(k/(n p_n))} \left( \frac{\widehat{\pi}_n}{\pi_n}-1 \right)\stackrel{d}{\longrightarrow} \Gamma. 
\]
\end{proposition}
\begin{remark}[Comparison with the heavy-tailed setting] In the heavy-tailed case, according to Section~5 in~\cite{daogirstu2018}, $\pi_n\approx p_n/(\gamma^{-1}-1)$. An estimator of $\pi_n$ is then $\widehat{\pi}_n=p_n/(\widehat{\gamma}_n^{-1}-1)$. In this setting, it is straightforward to obtain, under a suitable bias condition 
when
$\sqrt{k}(\widehat{\gamma}_n-\gamma)\to \Gamma$, 
that
\[
\sqrt{k} \left( \frac{\widehat{\pi}_n}{\pi_n}-1 \right)\stackrel{d}{\longrightarrow} \frac{\Gamma}{\gamma(1-\gamma)}. 
\]
The estimator $\widehat{\pi}_n$ therefore converges at a slightly faster rate in the heavy-tailed model. The slower 
speed
of convergence in the short-tailed framework is due to the presence of the quantity $\widehat{x}^{\star} - \widehat{q}_{1-p_n}^{\star}=-\widehat{a}(n/k) \widehat{\gamma}_n^{-1} (k/(n p_n))^{\widehat{\gamma}_n}$ in the numerator of $\widehat{\pi}_n$, whose rate of convergence to $x^{\star} - q_{1-p_n}$ is precisely $\log(k/(n p_n))/\sqrt{k}$, as obtained in Proposition~\ref{prop:estim_pin}.
\end{remark}
%

%
%
\section{Simulation study} \label{sec:fin}
%
%
The finite-sample performance of the 
proposed extreme expectile estimators
is illustrated here through a simulation study. 
Our setup first considers three models for i.i.d.~data: 
\begin{enumerate}[label=(\roman*), wide, labelindent=0pt]
\item The $X_t$ have a Beta distribution, whose density function is
\[
f(x|\alpha,\beta)=\frac{x^{\alpha-1}(1-x)^{\beta-1}}{\text{B}(\alpha,\beta)},\quad 0\leq x\leq 1. 
\]
Here $\text{B}(\alpha,\beta)$ is the Beta function and the shape parameters are set as $\alpha=3$ and $\beta=5/2$. 
In this model, the extreme value index is $\gamma=-2/5$ and the upper endpoint is $x^*=1$.
\item The $X_t$ have a short-tailed power-law distribution, whose distribution function is
\[
F(x|x^*,K,\alpha)= 1 - K(x^*-x)^\alpha, \quad x^*-K^{-1/\alpha}\leq x\leq x^*.
\]
Here $x^{\star}$, $K$ and $\alpha$ are the distribution endpoint, a positive constant and the shape parameter, respectively, which have been set as $x^*=5$, $K=1/3$ and $\alpha=3$, so that 
$\gamma=-1/3$. 
\item The $X_t$ have a GEV distribution, whose distribution function is
\[
F(x|\gamma)=\exp(-(1+\gamma x)^{-1/\gamma}), \quad 1+\gamma x >0.
\]
We set the extreme value index $\gamma=-1/3$, so that the upper endpoint is $x^*=-1/\gamma=3$.
\end{enumerate}
We then consider the following three time series models, in which $\Phi$ denotes the standard normal distribution function and $Y_t$ is the AR(1) process defined as $Y_{t+1}=\varrho Y_t + \sqrt{1-\varrho^2} \, \varepsilon_t$, with independent standard normal innovations $\varepsilon_t$, and where $\varrho\in (-1,1)$: 
\begin{enumerate}[label=(\roman*), wide, labelindent=0pt, resume]
\item $X_t=q_X(\Phi(Y_t))$, where $q_X$ is the quantile function corresponding to the Beta distribution defined in (i), and where the correlation parameter is $\varrho=0.95$.
\item $X_t=q_X(\Phi(Y_t))$, where $q_X$ is the quantile function corresponding to the short-tailed power-law distribution defined in (ii), and where the correlation parameter is $\varrho=0.5$.
\item $X_t=q_X(\Phi(Y_t))$, where $q_X$ is the quantile function corresponding to the GEV distribution defined in (iii), and where the correlation parameter is $\varrho=0.8$.
\end{enumerate}
The EVI and 
endpoints of models (iv)-(vi) are those of models (i)-(iii), respectively, and the time series models (iv)-(vi) are geometrically $\beta$-mixing (and in particular geometrically $\alpha-$mixing) since the linear AR(1) process $(Y_t)$ is so. We consider 
sample sizes $n=150, 300, 500$ and 
aim to predict expectiles of extreme level $\tau_n'=1-p_n=1-1/n= 0.9933, 0.9967, 0.9980$. 
As the true expectile values cannot be given in closed form, they have been computed by intensive Monte Carlo simulations and are reported in Table~\ref{tab:true_val_extreme}.
\begin{table}[b!]
\begin{center}
\begin{tabular}{cccc}
\toprule
Model & $\tau'_n=0.9933$ & $\tau'_n=0.9967$ & $\tau'_n=0.9980$\\
\midrule
(i), (iv) & $0.8571$ & $0.8814$ & $0.8968$ \\
(ii), (v) & $4.5284$ & $4.5939$ & $4.6372$ \\
(iii), (vi) & $1.9523$ & $2.1020$ & $2.2000$\\
\bottomrule
\end{tabular}
\caption{Values of the expectile $\xi_{\tau'_n}$ obtained through intensive Monte Carlo simulations for $\tau'_n=1-1/n$, with $n=150, 300, 500$.}
\label{tab:true_val_extreme}
\end{center}
\end{table}

We simulate $M=10{,}000$ samples of $n$ observations from each model and compare the purely empirical (LAWS) estimator $\widehat{\xi}_{\tau_n'}$ in \eqref{eq:empirical_LAWS}, the extrapolating LAWS estimators $\widehat{\xi}^{\star}_{\tau_n'}$ in \eqref{eq:extrapo_LAWS} obtained by setting $\widehat{\sigma}_n = \widehat{a}(1/\widehat{\overline{F}}_n(\widehat{\xi}_{\tau_n}))$, its alternative version $\overline{\xi}^{\star}_{\tau_n'}$ obtained with $\widehat{\sigma}_n = \widehat{a}((1-\tau_n)^{-1}) \times ((1-\tau_n)/\widehat{\overline{F}}_n(\widehat{\xi}_{\tau_n}))^{\widehat{\gamma}_n}$ 
in view of the approximation $a(1/\overline{F}(\xi_{\tau_n})) \approx ((1-\tau_n)/\overline{F}(\xi_{\tau_n}))^{\gamma} a((1-\tau_n)^{-1})$ that follows from the regular variation property of the scale function $a$, and the extrapolating 
quantile-based (QB)
estimator $\widetilde{\xi}^{\star}_{\tau_n'}$ in \eqref{eq:extrapo_QB}. In these last three estimators, $(\widehat{a}(n/k),\widehat{\gamma}_n)$ are either the pair of 
GP-pseudo-ML
estimators of $(a(n/k),\gamma)$~\citep[see][Section~3.4]{haafer2006} based on the top $k$ observations in the sample, or their versions based on the~\cite{dekeinhaa1989} Moment estimator~\citep[see][Section~3.5 and~4.2]{haafer2006}. 
We set throughout $\tau_n=1-k/n$, let the effective sample size $k$ range from $1\%$ up to $25\%$ of the total sample size $n$, and record Monte Carlo approximations of the relative bias, variance and Mean Squared Error (MSE) of the estimators as a function of $k$.

Results are reported in Figures 
B.1--B.6 in Section~B.1. In each figure the relative bias, variance and MSE are displayed from left to right, and results related to sample sizes $n=150,300,500$ are shown from top to bottom. For the sake of brevity we only report below in Figure~\ref{fig:simul:beta} 
the
results obtained with the Beta distribution, 
for the
sample size $n=300$ that we will also consider in our real data analysis, 
but we discuss the conclusions from the full set of models in 
Section~B.1. The Beta model corresponds to a case in which 
the
extreme value bias is present (unlike in the power-law setting, which is a transformation of a pure Pareto model) but not too disruptive in small samples (unlike in the case of the GEV distribution, which should be seen as difficult from that perspective). On the basis of the bias, the empirical estimator and extrapolating QB estimator tend to underestimate the true expectile along the entire range of 
$k$ values,
while the extrapolating LAWS estimator tends to overestimate the true expectile (at least when the scale and shape parameters are estimated via 
GP-pseudo-ML).
From the variance point of view, the extrapolating QB estimator is overall best among the estimators we consider, with the extrapolating LAWS estimators having large variance for small values of 
$k$.
Variability of the estimates seems to be highest when the data 
comes from time series. This conclusion carries over to the MSE: based on this criterion, the extrapolating QB estimator overall performs best, as expected from Remark~\ref{rmk:rates}, with the extrapolating LAWS estimator sometimes outperforming the extrapolating QB estimator for effective sample fractions larger than~$20\%$. In general, both extrapolating QB and LAWS estimators seem to perform remarkably well relative to the purely empirical expectile estimator, 
given the small sample size in this 
study.
%

%
%
\section{Application to forecast verification and comparison}\label{sec:data}
\begin{landscape}
\begin{figure}[h!] 
\vspace*{-1cm}
    \centering
        \includegraphics[width=6cm, height=5.5cm]{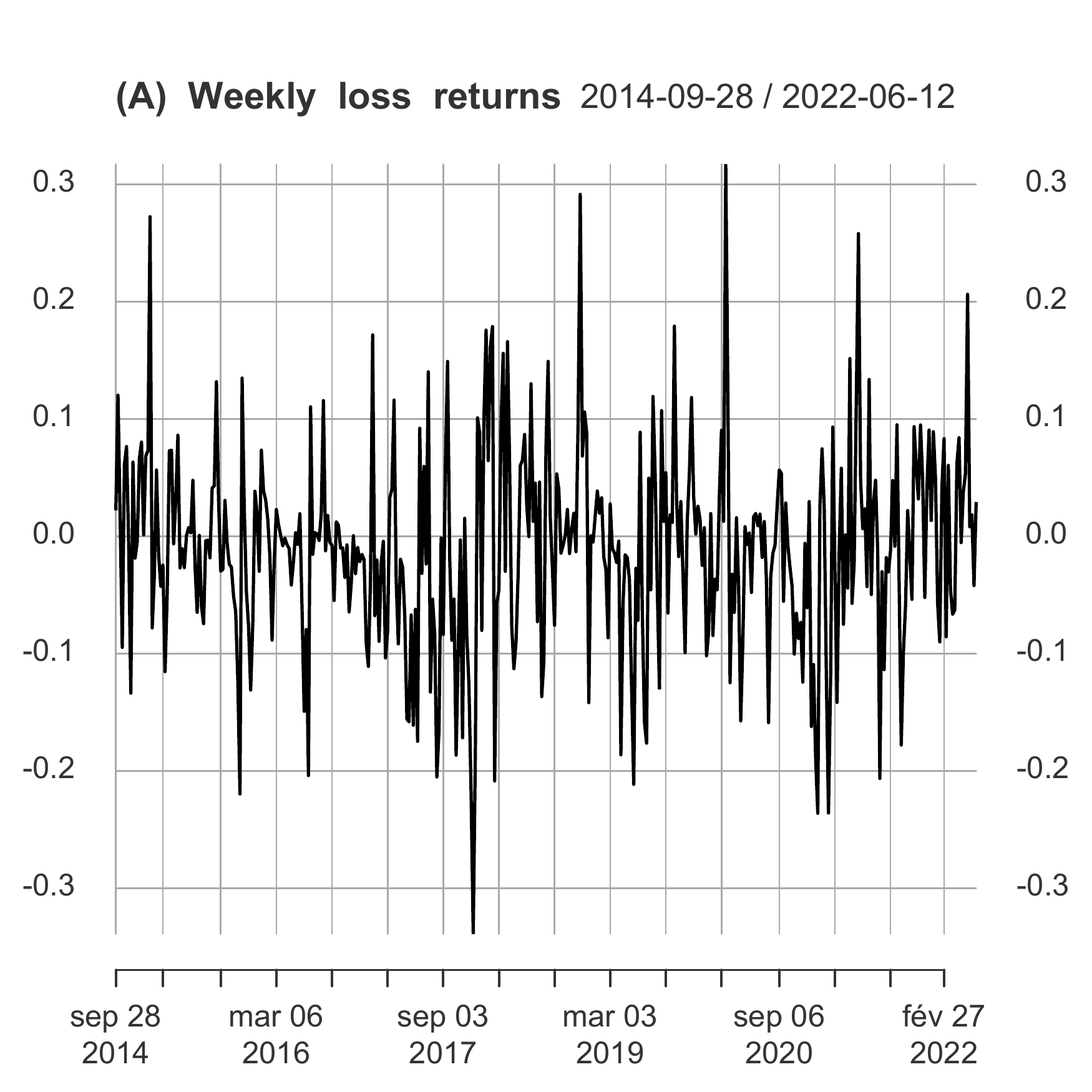}  
        \includegraphics[width=5.5cm, height=5cm]{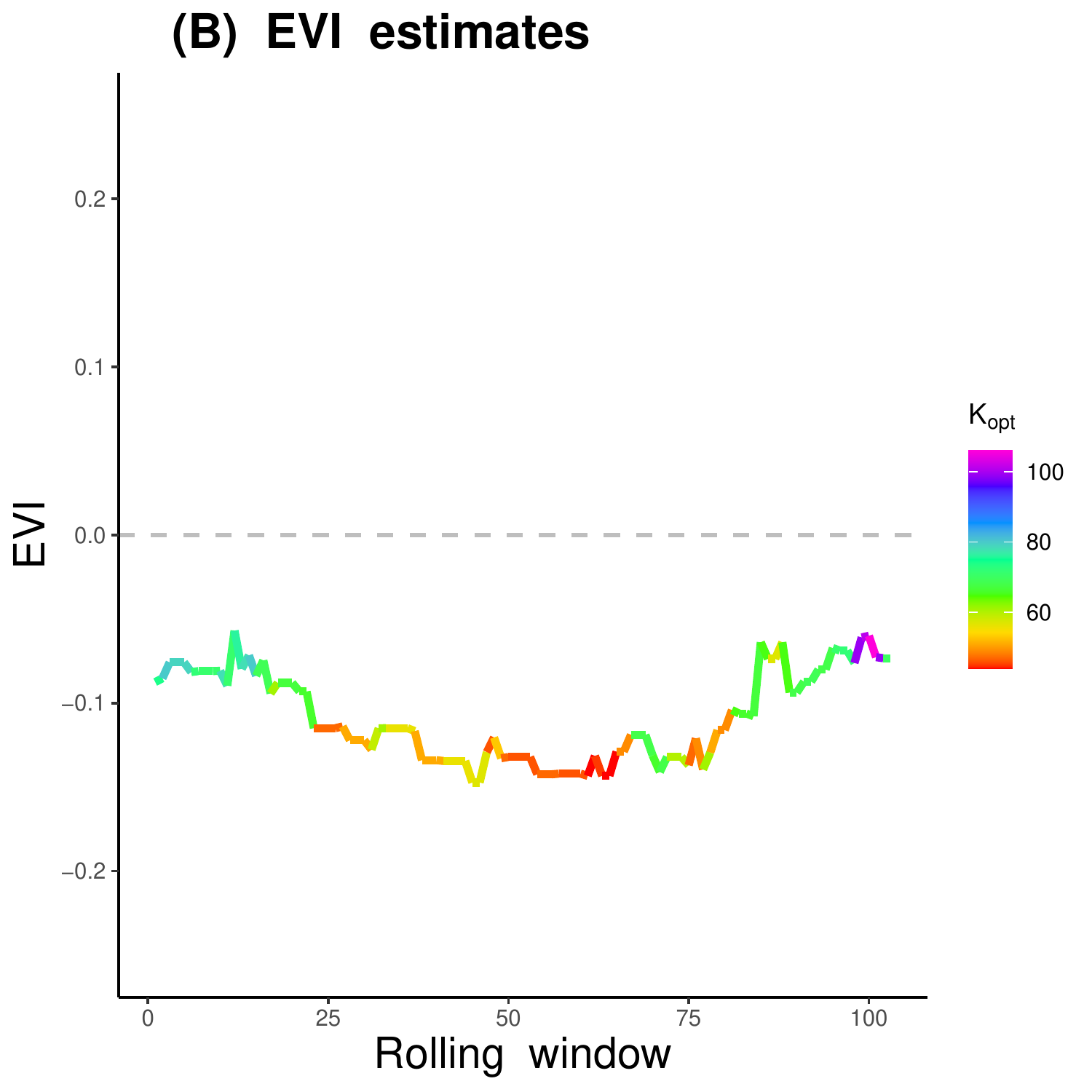} 
        \includegraphics[width=5.5cm, height=5cm]{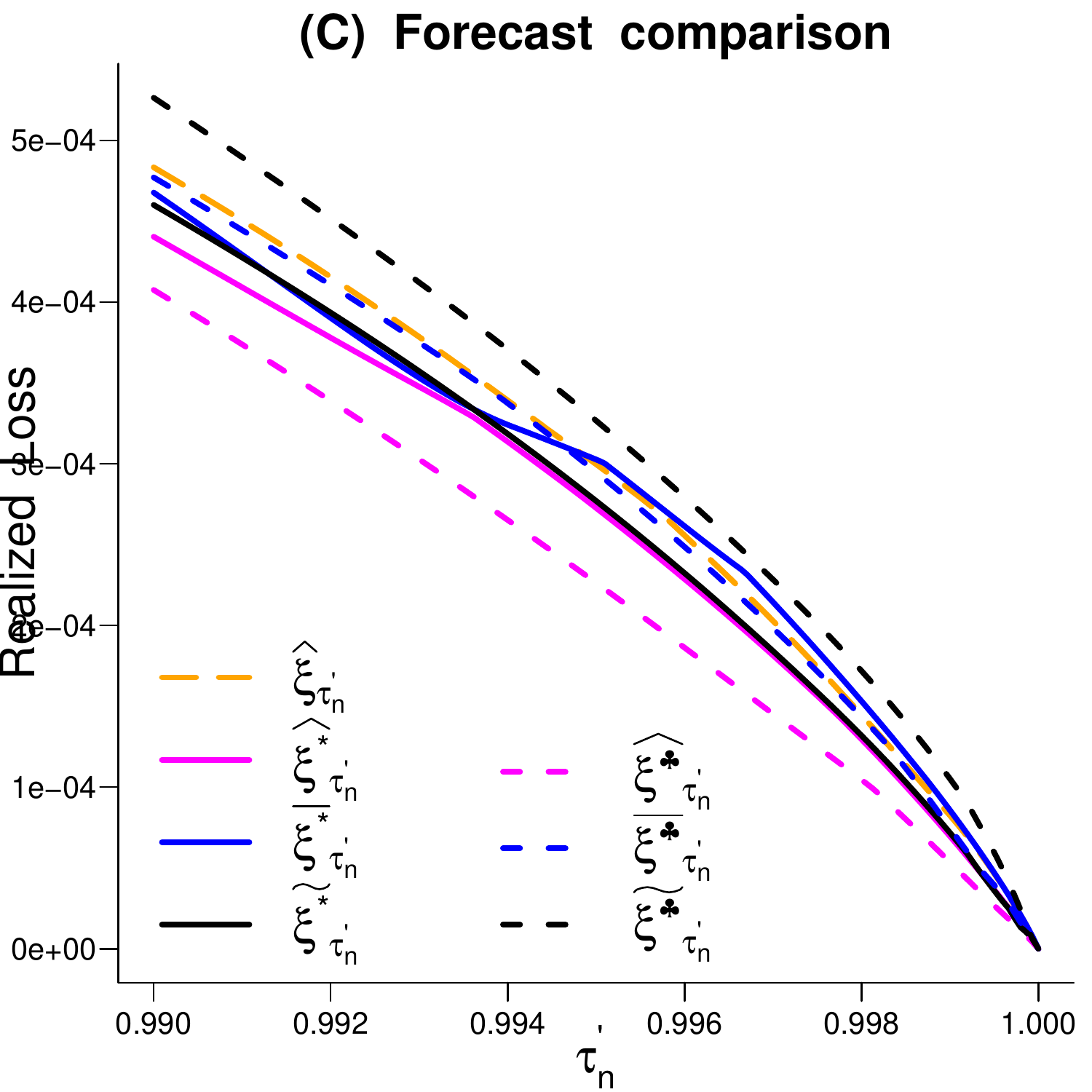} 
        \includegraphics[width=5.5cm, height=5cm]{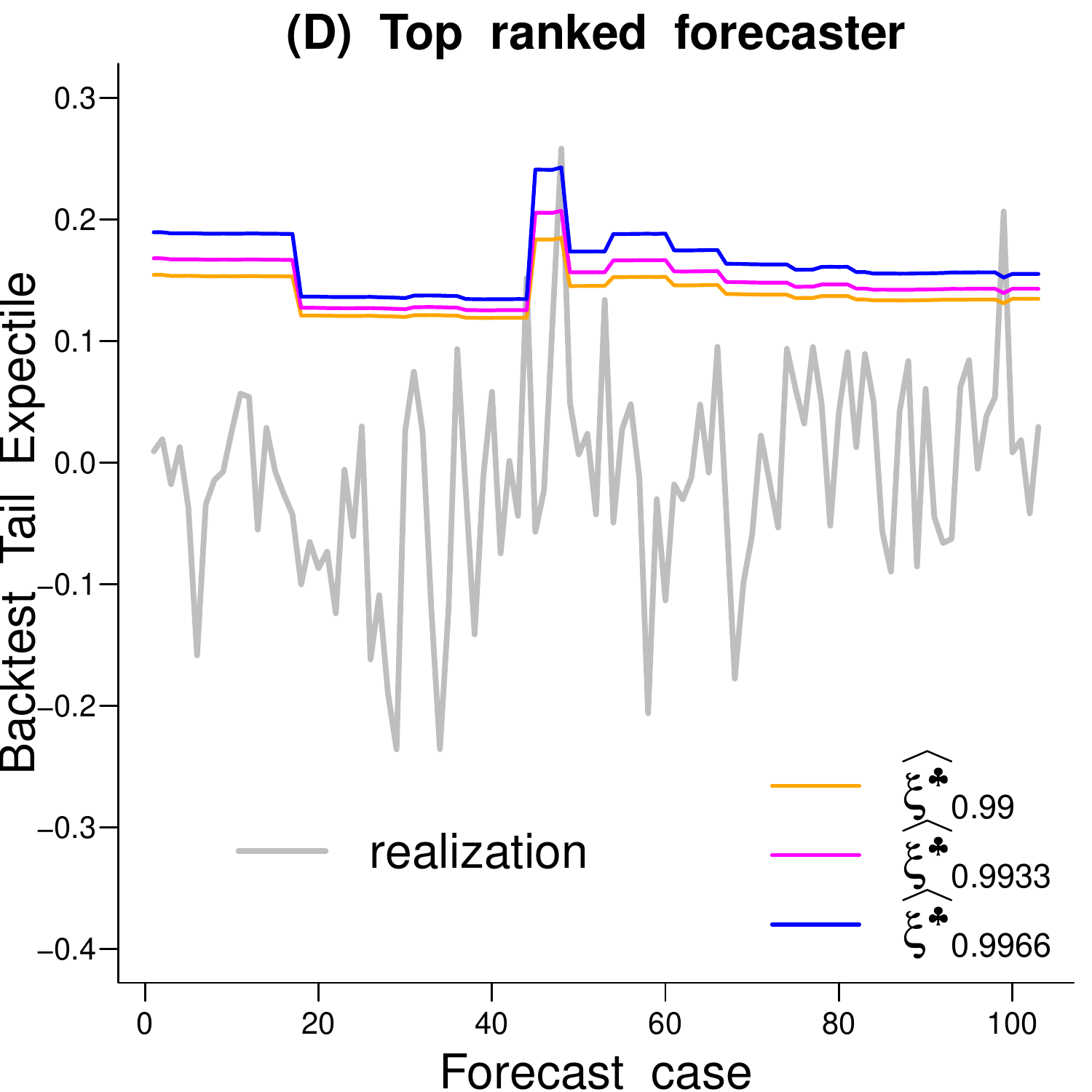} 
        \includegraphics[width=5.5cm, height=5cm]{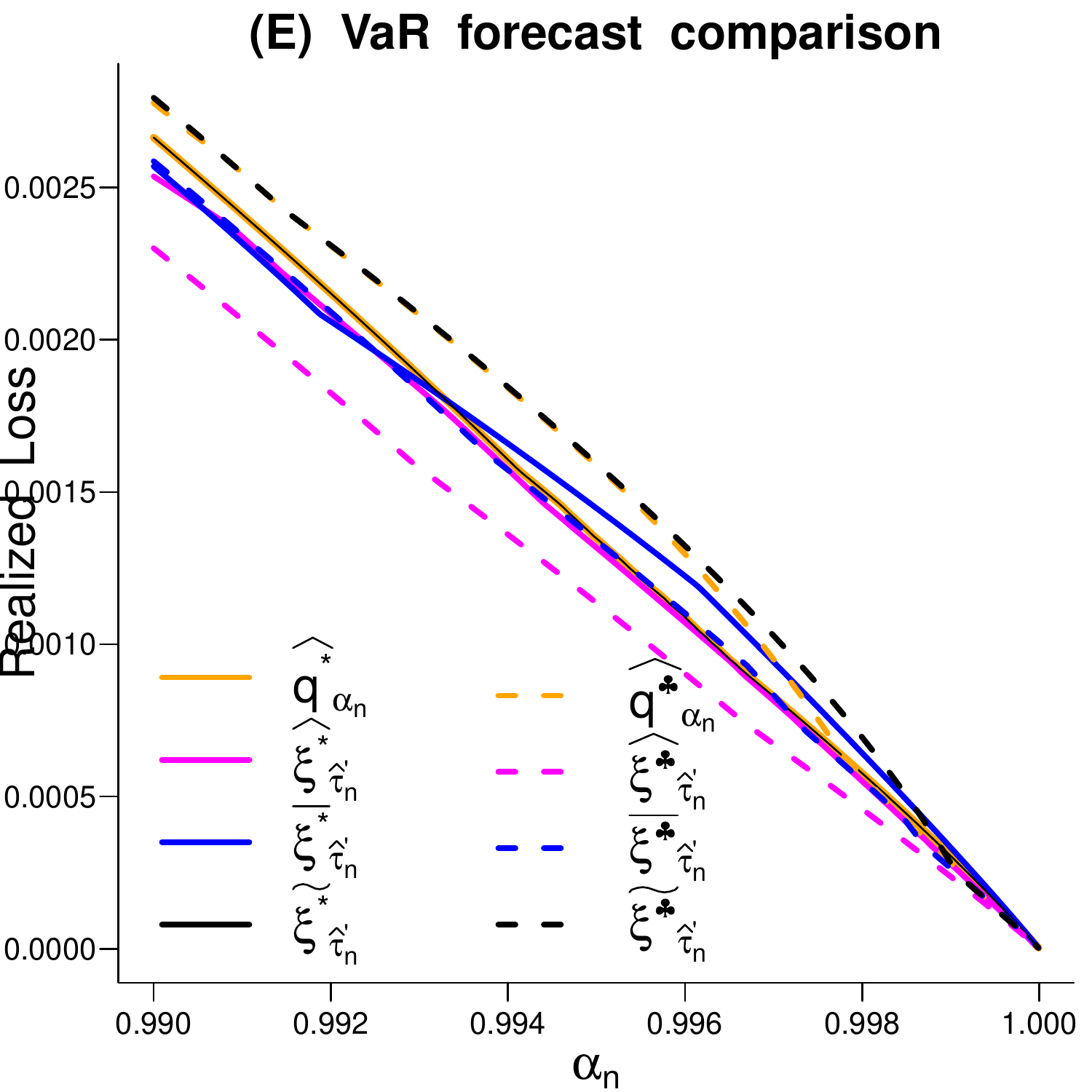} 
        \includegraphics[width=5.5cm, height=5cm]{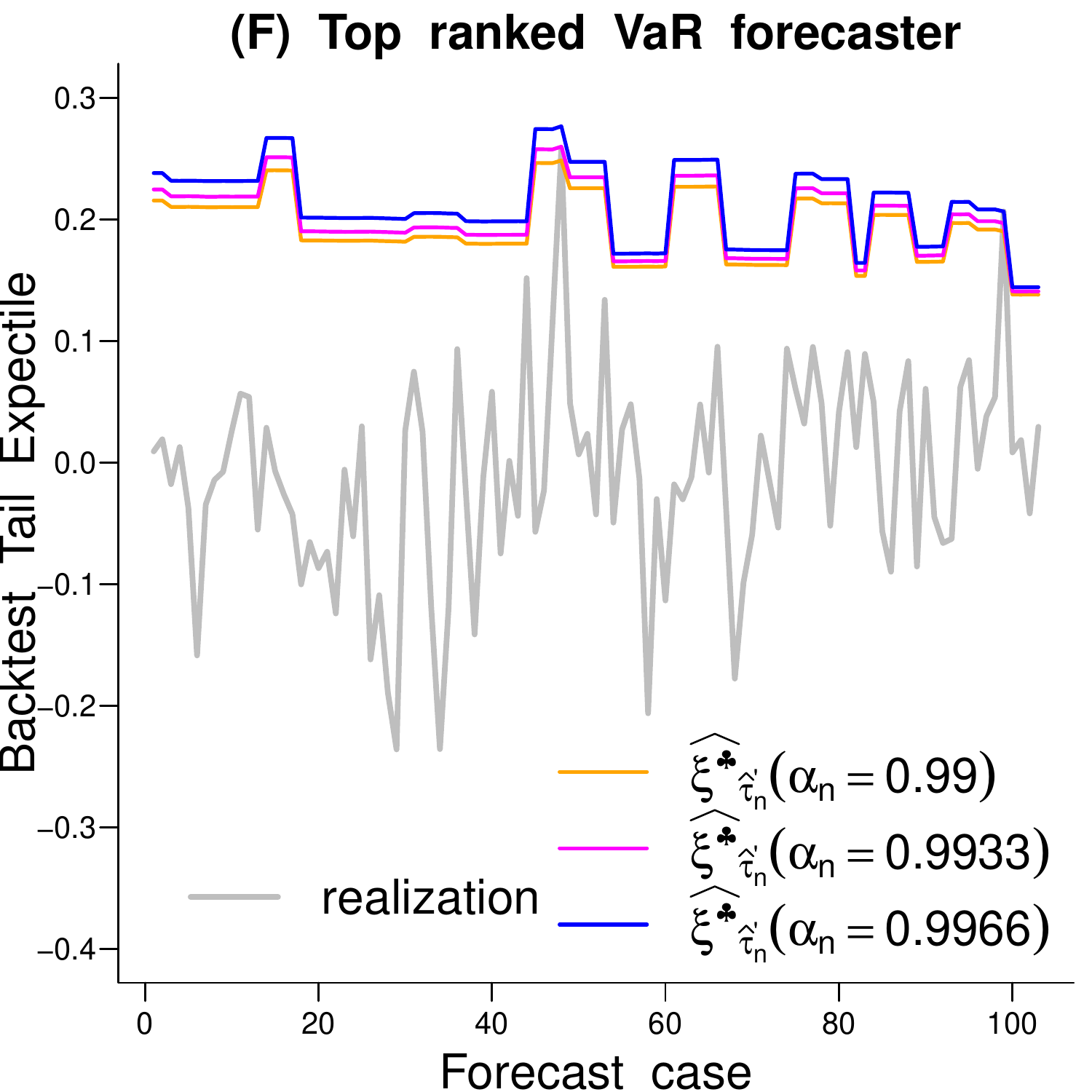} 

\caption{
(A)~Bitcoin weekly loss returns from September 28th, 2014, to June 12th, 2022. 
%
(B)~GP-pseudo-ML estimates of $\gamma$ over the 103 rolling windows. 
(C)~Realized loss function $\tau_n'\mapsto \overline{L}^{(m)}_{\tau_n'}$ evaluated at $\tau_n'\in [0.99, 1]$ 
for the seven competing forecasters   $\xi^{(1)}_t=\widehat\xi_{\tau_n'}$ (orange), $\xi^{(2)}_t =\widehat{\xi}_{\tau_n'}^{\star}$ (solid magenta),
  $\xi^{(3)}_t=\overline{\xi}_{\tau_n'}^{\star}$ (solid blue), $\xi^{(4)}_t =\widetilde{\xi}_{\tau_n'}^{\star}$ (solid black), $\xi^{(5)}_t =\widehat{\xi}_{\tau_n'}^{\clubsuit}$ (dashed magenta), $\xi^{(6)}_t =\overline{\xi}_{\tau_n'}^{\clubsuit}$ (dashed blue), and $\xi^{(7)}_t =\widetilde{\xi}_{\tau_n'}^{\clubsuit}$ (dashed black). 
(D)~The top-ranked optimal point forecaster 
$\widehat{\xi}_{\tau_n'}^{\clubsuit}$ of the risk measure $\xi_{\tau_n'}$, for $\tau_n'=0.99$ (orange), $\tau_n'=0.9933$ (magenta) and $\tau_n'=0.9966$ (blue), 
along with the realizing observations (gray).
(E)~Realized loss function $\alpha_n\mapsto \overline{L}^{(m)}_{\alpha_n}$ evaluated at $\alpha_n\in [0.99, 1]$ 
for the eight competing forecasters
$q^{(1)}_t:=\widehat{q}_{\alpha_n}^{\star}$ (solid orange), 
$q^{(2)}_t :=\widehat{\xi}_{\widehat\tau_n'}^{\star}$ (solid magenta), $q^{(3)}_t:=\overline{\xi}_{\widehat\tau_n'}^{\star}$ (solid blue), 
$q^{(4)}_t :=\widetilde{\xi}_{\widehat\tau_n'}^{\star}$ (solid black), 
$q^{(5)}_t:=\widehat{q}_{\alpha_n}^{\clubsuit}$ (dashed orange),
$q^{(6)}_t :=\widehat{\xi}_{\widehat\tau_n'}^{\clubsuit}$ (dashed magenta),
$q^{(7)}_t :=\overline{\xi}_{\widehat\tau_n'}^{\clubsuit}$ (dashed blue), and
$q^{(8)}_t :=\widetilde{\xi}_{\widehat\tau_n'}^{\clubsuit}$ (dashed black).
(F) The top-ranked optimal point forecaster
$\widehat{\xi}_{\widehat\tau_n'}^{\clubsuit}$ of the risk measure $q_{\alpha_n}$, for $\alpha_n=0.99$ (orange), $\alpha_n=0.9933$ (magenta) and $\alpha_n=0.9966$ (blue), 
along with the realizing observations (gray).}
\label{fig:BTC:n300}
\end{figure}
\end{landscape}

%
%
In this section, we apply our LAWS and QB 
methods to estimate tail expectile risk for Bitcoin (BTC-USD), a peer-to-peer digital decentralized cryptocurrency. At the end of September 2014, Bitcoin had volatility seven times greater than gold, eight times greater than the S\&P~500, and 18 times greater than the US dollar. Although the growth of Bitcoin prices has been often described as an economic bubble, the COVID-19 crisis has sparked substantial investment in this digital currency as an alternative to conventional asset classes. Similarly to the 14 major companies and financial institutions that we have explored in Figure~\ref{fig:data:EVI}, we will provide evidence that short-tailed returns 
occur over (relatively) short time periods for this 
cryptocurrency as well. To do so and to assess its associated extreme financial risk accordingly,
we construct a time series of weekly 
negative log-returns
from averaged daily Bitcoin closing prices within the corresponding week, from September 28th, 2014, to June 12th, 2022. 
The time series of loss returns is represented in 
Figure~\ref{fig:BTC:n300}~(A).

We consider risk assessment from a forecasting perspective. With our knowledge of this week, the goal is to give the best possible point estimate of the expectile risk measure $\xi_{\tau_n'}$ for the next week based on rolling windows of length $n = 300$. This 
results in 103 samples of size $n$ over the observed timeframe. For each sample $(X_1,\ldots,X_n)$, the EVI of the underlying distribution was estimated by means of the ML method for peaks over a high threshold $X_{n-k,n}$. The plot of the estimates obtained over the successive 103 rolling windows is given in Figure~\ref{fig:BTC:n300}~(B), where an appropriate $k$ is chosen, for each sample, by regarding the path of the 
$\gamma$ estimator
as a function of~$k$ and selecting the~$k$ value which corresponds to the median estimate over the most stable region of the path (this can be achieved by using 
the algorithm developed by~\cite{elmstu2017}).
This selection is highlighted in Figure~\ref{fig:BTC:n300}~(B) by a colour scheme, ranging from dark red (low) to dark violet (high). The final 
$\gamma$ estimates are found to be all negative in $[-0.147,-0.057]$. 
It should be noted that we have comfortably concluded the stationarity of the time series samples across all $T=103$ rolling windows, from our exploratory analysis.

Expectiles have recently received growing attention in quantitative risk management not only for their coherence as a tail risk measure, but also for their 
elicitability that corresponds to the existence of a natural methodology for 
forecast verification. According to Gneiting (2011) and Ziegel (2016), 
letting the random variable $X$ model the future observation of interest, 
$\xi_{\tau_n'}$ equals the optimal point forecast for $X$ given by the Bayes rule 
$
\xi_{\tau_n'}= \argmin_{\xi\in\mathbb{R}} \mathbb{E}\left[L_{\tau_n'}(\xi,X)\right]
$, 
under the  asymmetric quadratic scoring function
$$
L_{\tau_n'}: \, \mathbb{R}^2 \longrightarrow [0,\infty), \quad (\xi,x)\mapsto \eta_{\tau_n'}(x-\xi),
$$
where $L_{\tau_n'}(\xi,x)$ represents the loss or penalty when the point forecast $\xi$ is issued and the realization $x$ of $X$ materializes.
Following the ideas of Gneiting (2011) and Ziegel (2016),
the competing 
estimation procedures for $\xi_{\tau_n'}$ can be compared by using the scoring function~$L_{\tau_n'}$:
Suppose that, in $T$ forecast cases, we have point forecasts $\left(\xi^{(m)}_1,\ldots,\xi^{(m)}_T\right)$ and realizing observations $(x_1,\ldots,x_T)$, where the index $m$ numbers the competing forecasters that are computed at each forecast case $t=1,\ldots,T$. In the assessment, we compare the purely empirical expectile $\xi^{(1)}_t:=\widehat\xi_{\tau_n'}$ in~\eqref{eq:empirical_LAWS} with the direct extrapolating LAWS estimator $\xi^{(2)}_t :=\widehat{\xi}_{\tau_n'}^{\star}$ in~\eqref{eq:extrapo_LAWS} and its alternative version $\xi^{(3)}_t:=\overline{\xi}_{\tau_n'}^{\star}$ described in Section~\ref{sec:fin}, and with the indirect 
QB
extrapolating estimator $\xi^{(4)}_t :=\widetilde{\xi}_{\tau_n'}^{\star}$ in~\eqref{eq:extrapo_QB}, all of them being based on the 
GP-pseudo-ML
estimators $(\widehat{a}(n/k),\widehat{\gamma}_n)$ of $(a(n/k),\gamma)$. When the 
Moment
estimators $(\widehat{a}(n/k),\widehat{\gamma}_n)$ are used instead of the ML estimators, the corresponding three extrapolating forecasters $\widehat{\xi}_{\tau_n'}^{\star}$, $\overline{\xi}_{\tau_n'}^{\star}$ and $\widetilde{\xi}_{\tau_n'}^{\star}$ will be denoted in the sequel by replacing ``$\star$'' with ``{\scriptsize $\clubsuit$}'' to define
$\xi^{(5)}_t :=\widehat{\xi}_{\tau_n'}^{\clubsuit}$,
$\xi^{(6)}_t :=\overline{\xi}_{\tau_n'}^{\clubsuit}$ and
$\xi^{(7)}_t :=\widetilde{\xi}_{\tau_n'}^{\clubsuit}$.
The seven competing point estimates can then be ranked in terms of their average scores (the lower the better):
\begin{equation*}\label{barL}
\overline{L}^{(m)}_{\tau_n'} = \frac{1}{T}\sum_{t=1}^T  L_{\tau_n'}\left(\xi^{(m)}_t,x_t\right), \quad m=1,\ldots,7.
\end{equation*}
The computation of the different extrapolated expectile estimators requires, like the EVI estimators, the determination of the optimal value of the effective sample size $k$. 
By balancing the potential estimation bias and variance, a usual practice in extreme value theory is to choose $k$ from the first stable region of the plots [see, {\it e.g.}, Section~3 in de Haan and Ferreira (2006)]. This is achieved by using the path stability procedure for~$\gamma$ estimation. However, to achieve optimal point forecasts $\left(\xi^{(m)}_1,\ldots,\xi^{(m)}_T\right)$ for the future observation~$X$, this requires the use of~$k$ values that minimize their associated realized loss $\overline{L}^{(m)}_{\tau_n'}\equiv \overline{L}^{(m)}_{\tau_n'}(k)$, for m=2,\ldots,7. 
Doing so, we obtain the final values of $\overline{L}^{(m)}_{\tau_n'}$ graphed in Figure~\ref{fig:BTC:n300}~(C), as functions of the extreme level $\tau_n'\in [0.99, 1]$, for the seven competing estimators. It can be seen that the 
LAWS-Moment
estimator $\widehat{\xi}_{\tau_n'}^{\clubsuit}$ (dashed magenta) is the best forecaster uniformly in $\tau_n'$, followed by the LAWS-ML estimator $\widehat{\xi}_{\tau_n'}^{\star}$ (solid magenta) and then the QB-ML estimator $\widetilde{\xi}_{\tau_n'}^{\star}$ (solid black). The remaining three extrapolating estimators do not seem, for this particular choice of $T=103$ rolling windows of length $n=300$, to outperform the naive sample expectile $\widehat\xi_{\tau_n'}$ (dashed orange). The values of the top-ranked forecaster $\widehat{\xi}_{\tau_n'}^{\clubsuit}$, computed on the 103 successive rolling windows for the extreme levels $\tau_n'\in\{0.99, 0.9933, 0.9966\}$, are displayed in Figure~\ref{fig:BTC:n300}~(D), along with the realizing observation at each forecast case. The point forecasts seem to smoothly increase with $\tau_n'$ 
approaching the worst expected (finite) losses at $\tau_n' = 1$.
From the perspective of pessimistic decision making, 
the forecasts obtained at the lower level $\tau_n' =0.99$ (orange curve) are already cautious since they do lie almost overall beyond the range of the data: This is mainly due to the short-tailed nature of Bitcoin data that is closer to light-tailedness.

Extreme expectiles can also serve as a useful tool for estimating the conventional Value at Risk (VaR) itself. Stated differently, if the statistician or the practitioner wishes to forecast a coherent expectile $\xi_{\tau_n'}$ that has the same probabilistic interpretation as an extreme quantile $q_{\alpha_n}$, for a pre-specified tail probability level $\alpha_n$, a natural way of doing so is to select the asymmetry level $\tau_n'$ so that $\xi_{\tau_n'}\equiv q_{\alpha_n}$. As justified in Section~\ref{sec:main:select}, such a $\tau_n'$ can be estimated by $\widehat\tau_n'= 1-\widehat{\pi}_n$. When substituting this estimated value in place of $\tau_n'$ in our $\xi_{\tau_n'}$ extrapolated estimators, the latter estimate the VaR $q_{\alpha_n}$ itself and can then be compared with the popular GP fit $\widehat{q}_{\alpha_n}^{\star}$ defined in~\eqref{eq:quantile_endpoint}. Here also, forecast verification and comparison is possible thanks to the elicitability property of quantiles (see {\it e.g.} Gneiting (2011)). Given that it is the quantile level $\alpha_n$ which is fixed in advance, the accuracy of the associated VaR forecasts is to be assessed by means of the realized loss
\begin{equation*}
\overline{L}^{(m)}_{\alpha_n} = \frac{1}{T}\sum_{t=1}^T  L_{\alpha_n}\left(q^{(m)}_t,x_t\right), \quad m=1,\ldots,8,
\end{equation*}
under the asymmetric piecewise linear scoring function
$$
L_{\alpha_n}: \, \mathbb{R}^2 \longrightarrow [0,\infty), \quad (q,x)\mapsto \varrho_{\alpha_n}(x-q),
$$
for the competing ML-based forecasters 
$q^{(1)}_t:=\widehat{q}_{\alpha_n}^{\star}$, 
$q^{(2)}_t :=\widehat{\xi}_{\widehat\tau_n'}^{\star}$, $q^{(3)}_t:=\overline{\xi}_{\widehat\tau_n'}^{\star}$, 
$q^{(4)}_t :=\widetilde{\xi}_{\widehat\tau_n'}^{\star}$, 
and their 
Moment-based
versions 
$q^{(5)}_t:=\widehat{q}_{\alpha_n}^{\clubsuit}$,
$q^{(6)}_t :=\widehat{\xi}_{\widehat\tau_n'}^{\clubsuit}$,
$q^{(7)}_t :=\overline{\xi}_{\widehat\tau_n'}^{\clubsuit}$, and
$q^{(8)}_t :=\widetilde{\xi}_{\widehat\tau_n'}^{\clubsuit}$.
The resulting realized losses $\overline{L}^{(m)}_{\alpha_n}$ are graphed in Figure~\ref{fig:BTC:n300}~(E), as functions of the quantile level $\alpha_n\in [0.99, 1]$, for the eight competing estimators of $q_{\alpha_n}\equiv\xi_{\tau_n'}$. It is remarkable that the best forecaster is still the 
LAWS-Moment
estimator $\widehat{\xi}_{\widehat\tau_n'}^{\clubsuit}$ (dashed magenta), followed by the LAWS-ML estimator $\widehat{\xi}_{\widehat\tau_n'}^{\star}$ (solid magenta). Most importantly, these expectile-based forecastors clearly outperform the usual GP-ML fit $\widehat{q}_{\alpha_n}^{\star}$ (solid orange) and GP-Moment fit $\widehat{q}_{\alpha_n}^{\clubsuit}$ (dashed orange), which is good news to practitioners whose concern is to assess the accuracy of forecasts. Figure~\ref{fig:BTC:n300}~(F) contrasts the evolution of the optimal point forecasts $\widehat{\xi}_{\widehat\tau_n'}^{\clubsuit}$, for the risk measure $q_{\alpha_n}$ at the extreme levels $\alpha_n\in\{0.99, 0.9933, 0.9966\}$, with the realizing 
observations.
By comparing these $\alpha_n$th quantile estimates with their expectile analogs from Figure~\ref{fig:BTC:n300}~(D) at the same asymmetry levels ($\alpha_n=\tau_n'$), it may be seen that expectiles are ultimately less conservative than quantiles, which empirically corroborates the theoretical result for short-tailed data in Proposition~2.2 by~\cite{beldib2017}. This more liberal expectile assessment of tail risk is indeed a consequence of the diversification principle satisfied by expectiles. 
Interestingly, the conservative LAWS-Moment (expectile-based) forecasts $\widehat{\xi}_{\widehat\tau_n'}^{\clubsuit}$, for $q_{\alpha_n}$ in Figure~\ref{fig:BTC:n300}~(F), seem also to be more sensitive to the variability of weekly losses compared with their analog forecasts $\widehat{\xi}_{\alpha_n}^{\clubsuit}$ for $\xi_{\alpha_n}$ in Figure~\ref{fig:BTC:n300}~(D).

%
\begin{figure}[h!]
\centering
%
%
\includegraphics[width=0.32\textwidth, page=4]{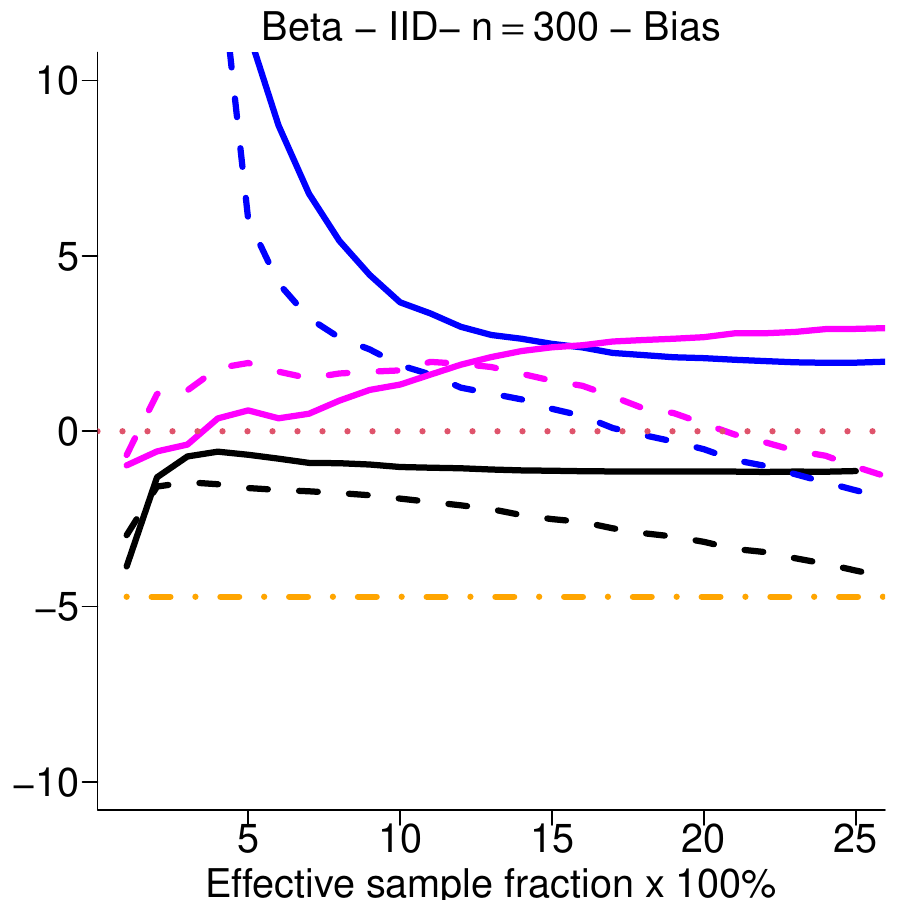}
\includegraphics[width=0.32\textwidth, page=5]{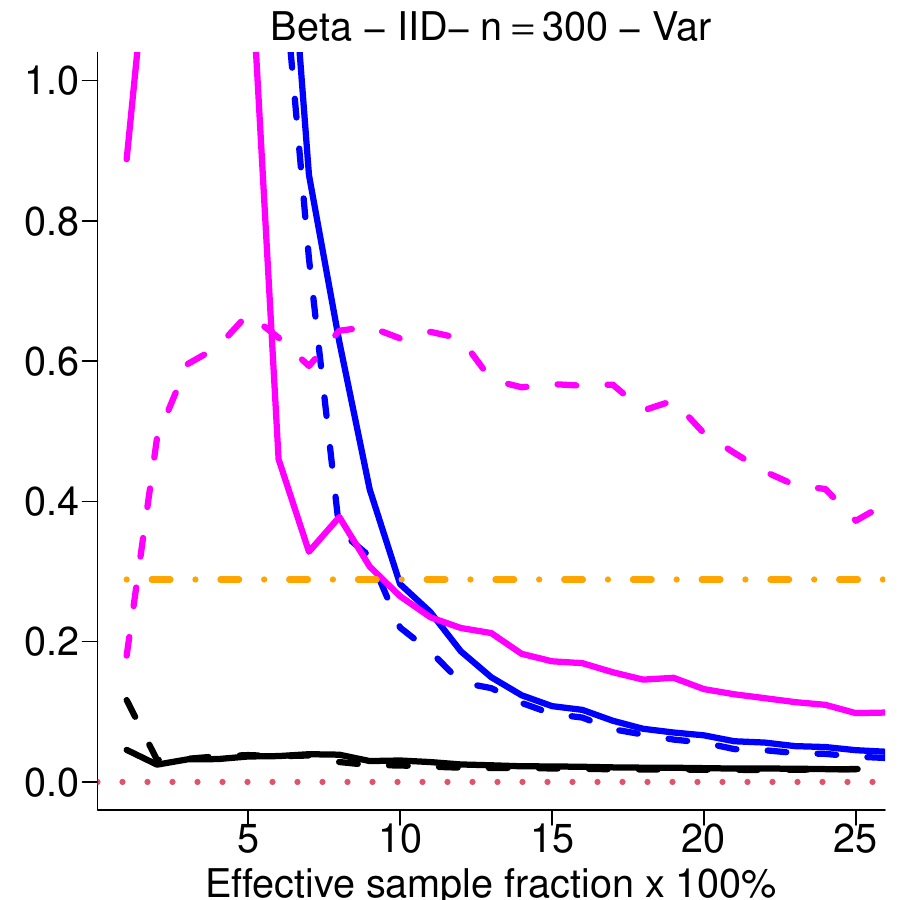}
\includegraphics[width=0.32\textwidth, page=6]{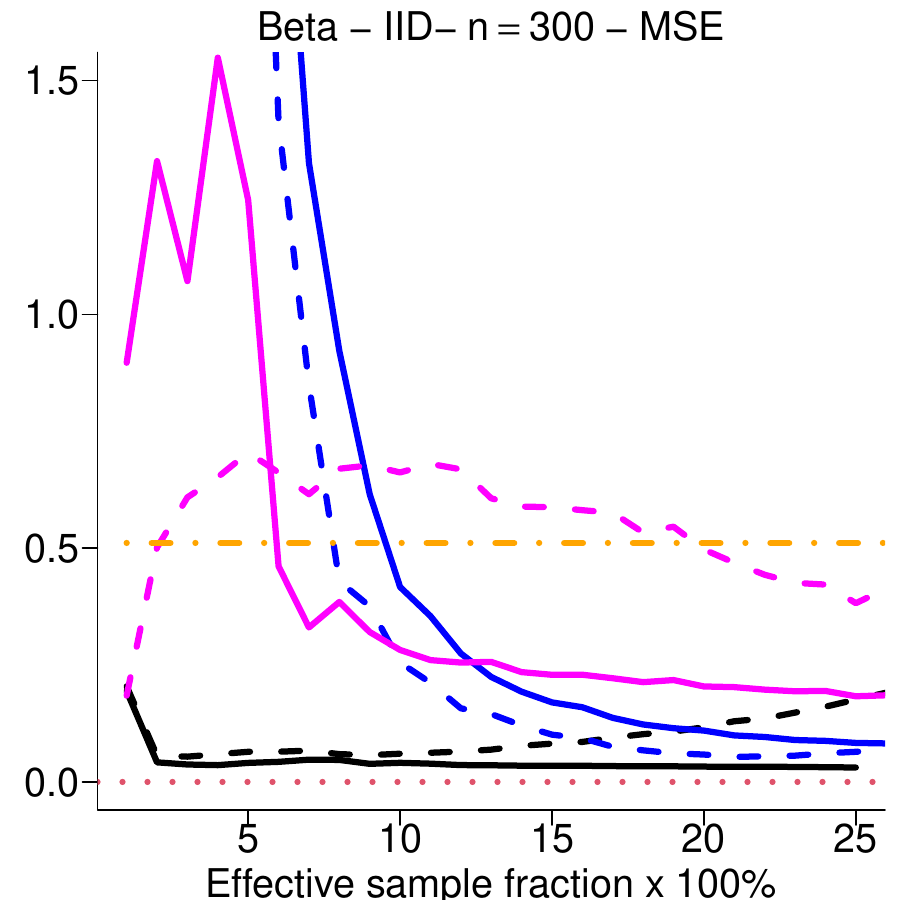} \\[10pt]
%
%
\includegraphics[width=0.32\textwidth, page=4]{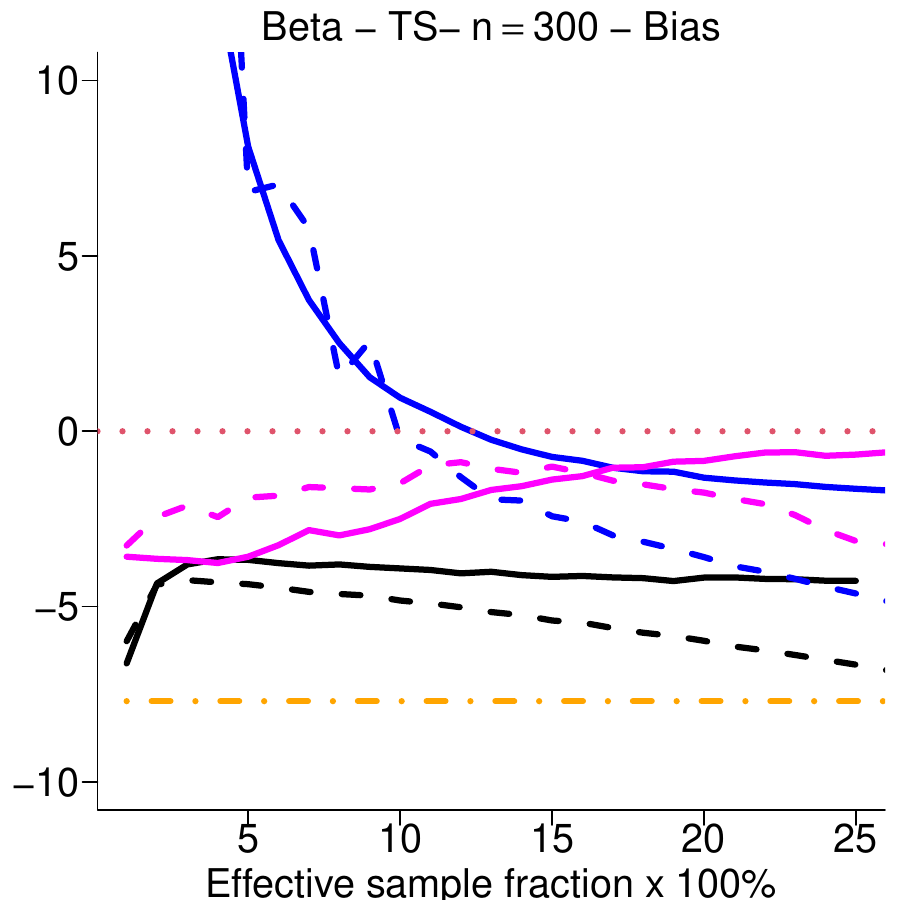}
\includegraphics[width=0.32\textwidth, page=5]{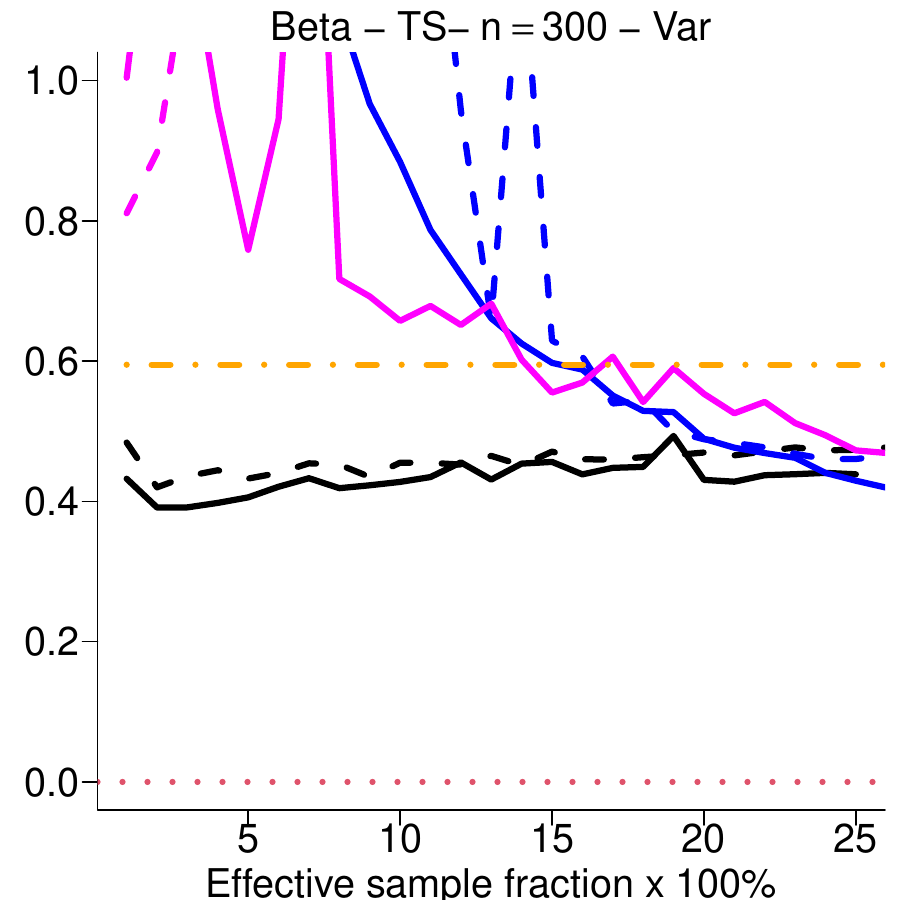}
\includegraphics[width=0.32\textwidth, page=6]{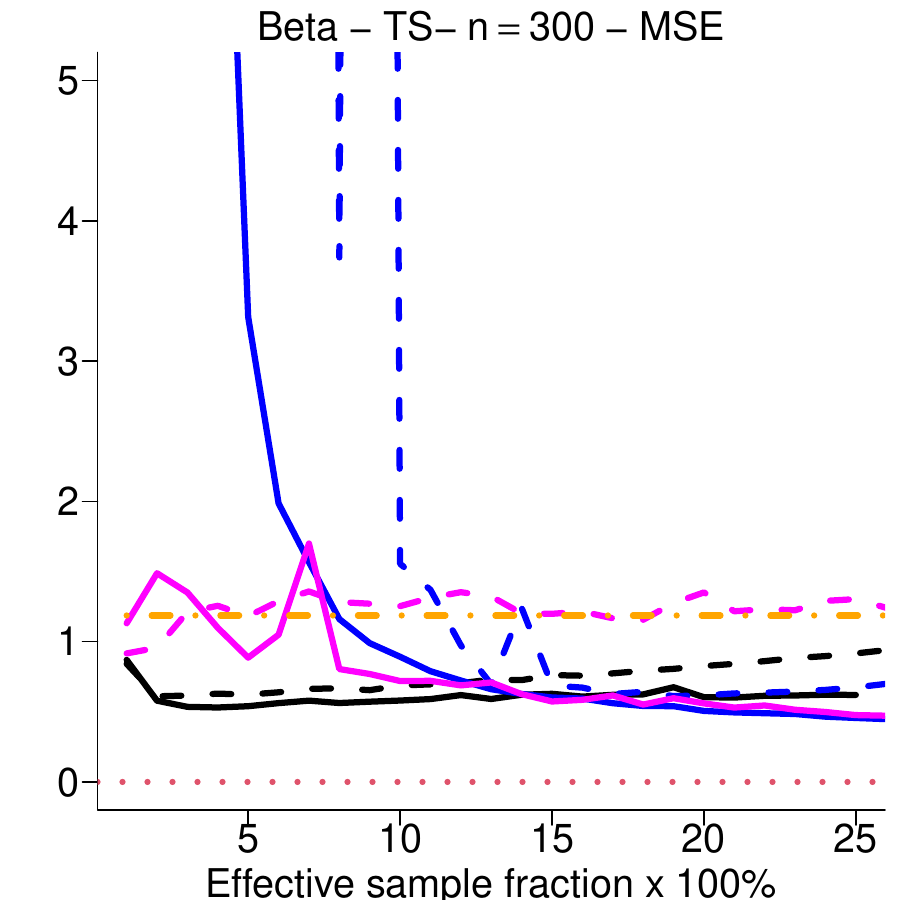} \\[10pt]
\caption{Empirical relative bias, variance and MSE (left, middle and right), multiplied by $100$, for the 
 $\xi_{\tau_n'}$ estimators
obtained with observations from a Beta distribution, $\tau_n'=1-1/n$ and 
$n=300$.  Empirical estimator $\widehat{\xi}_{\tau'_n}$ (orange line), extrapolating LAWS estimators $\widehat{\xi}^{\star}_{\tau'_n}$ (magenta lines) and $\overline{\xi}^{\star}_{\tau'_n}$ (blue lines), and extrapolating QB estimators $\widetilde{\xi}^{\star}_{\tau'_n}$ (black lines). The 
 extrapolating estimators based on the scale and shape parameter estimates from the GP-pseudo-ML are referred to using 
solid lines, and those based on the Moment estimators are referred to using dashed lines. Top: i.i.d.~data, bottom: nonlinear AR(1) data.}
\label{fig:simul:beta}
\end{figure}
%


%

\section*{Acknowledgments}

This research was supported by the French National Research Agency under the grants ANR-19-CE40-0013 and ANR-17-EURE-0010. S.A. Padoan is supported by the Bocconi Institute for Data Science and Analytics (BIDSA), Italy. A.~Daouia and G.~Stupfler acknowledge financial support from the TSE-HEC ACPR Chair and from an AXA Research Fund Award on ``Mitigating risk in the wake of the COVID-19 pandemic''.

\bibliographystyle{apalike}
\bibliography{biblio}

\end{document}